\documentclass[11pt,amsfonts]{article}
\usepackage{graphicx}
\usepackage{latexsym}
\usepackage{amssymb}
\usepackage{amsmath}
\usepackage{enumerate}
\usepackage{color}
\usepackage{layout}
\usepackage{eufrak}
\newtheorem{prop}{Proposition}
\newtheorem{lemma}{Lemma}

\newtheorem{corollary}{Corollary}
\newtheorem{theorem}{Theorem}
\newtheorem{remark}{Remark}

\def\real{{\mathord{{\rm I\kern-2.8pt R}}}}        % Fake blackboard bold R.
\def\inte{{\mathord{{\rm I\kern-2.8pt N}}}}

\def\sZZ{{\rm Z\kern-2.8ptem{}Z}}

\def\z{{\mathchoice
  {\sZZ}
  {\sZZ}
  {\rm Z\kern-0.30em{}Z}
  {\rm Z\kern-0.25em{}Z} }}
\def\sQQ{{\kern 0.27em \vrule height1.45ex width0.03em depth0em
          \kern-0.30em \rm Q}}
\def\qu{{\mathchoice
    {\sQQ}
    {\sQQ}
  {\kern 0.225em \vrule height1.05ex width0.025em depth0em \kern-0.25em \rm Q}
  {\kern 0.180em \vrule height0.78ex width0.020em depth0em \kern-0.20em \rm Q}
        }}
\def\sCC{{\kern 0.27em \vrule height1.45ex width0.03em depth0em
          \kern-0.30em \rm C}}
\def\complex{{\mathchoice
    {\sCC}
    {\sCC}
  {\kern 0.225em \vrule height1.05ex width0.025em depth0em \kern-0.25em \rm C}
  {\kern 0.180em \vrule height0.78ex width0.020em depth0em \kern-0.20em \rm C}
        }}

\newcommand{\ba}{\begin{array}}
\newcommand{\ea}{\end{array}}
\newcommand{\be}{\begin{equation}}
\newcommand{\ee}{\end{equation}}
\newcommand{\bea}{\begin{eqnarray}}
\newcommand{\eea}{\end{eqnarray}}
\newcommand{\beaa}{\begin{eqnarray*}}
\newcommand{\eeaa}{\end{eqnarray*}}

\def\b{\beta}

\def\z{\zeta}

\font\tenmath=msbm10 \font\sevenmath=msbm7 \font\fivemath=msbm5
\newfam\mathfam \textfont\mathfam=\tenmath
\scriptfont\mathfam=\sevenmath \scriptscriptfont\mathfam=\fivemath

\def \b{\noindent}

\def \={{\buildrel {\rm (law)} \over =}}

\def\qed{ \hfill \vrule width.25cm height.25cm depth0cm\smallskip}

\newcommand{\basa}{\begin{assumption}}
\newcommand{\easa}{\end{assumption}}

\newcommand{\bas}{\begin{assum}}
\newcommand{\eas}{\end{assum}}

\def\Re{\hbox{\rm Re$\,$}}

\newcommand{\ignore}[1]{}
\begin{document}

\renewcommand{\thefootnote}{\fnsymbol{footnote}}

\renewcommand{\thefootnote}{\fnsymbol{footnote}}

\title{Multidimensional Selberg theorem and fluctuations of the zeta zeros via Malliavin calculus }

\author{Ciprian A. Tudor $^{1}$\\, 
%\footnote{Supported by the CNCS grant PN-II-ID-PCCE-2011-2-0015.  }\vspace*{0.1in} \\
 $^{1}$ Laboratoire Paul Painlev\'e, Universit\'e de Lille 1\\
 F-59655 Villeneuve d'Ascq, France.\\
%Academy of Economical Studies, Bucharest, Romania. \vspace*{0.1in}\\
 \quad tudor@math.univ-lille1.fr\vspace*{0.1in}}
\maketitle

\begin{abstract}
We give new contributions on the distribution of the zeros of the Riemann zeta function by using the techniques of the Malliavin calculus. In particular, we obtain the error bound in the multidimensional  Selberg' s central limit theorem concerning the zeta zeros on the critical line and we discuss some consequences concerning the asymptotic behavior of the mesoscopic  fluctuations of the zeta zeros. 
\end{abstract}

\vskip0.2cm

{\bf 2010 AMS Classification Numbers:}  60F05, 60H05, 11M06.

 \vskip0.2cm


\section{On the Riemann zeta function and Selberg's theorem}

The Riemann zeta function is usually defined, for $\Re s>1$, as
\begin{equation}
\label{zeta1}
\zeta (s) =\sum_{n\geq 1} \frac{1}{ n^{s}}
\end{equation} 
and for $\Re s \leq 1$, as an analytic continuation of (\ref{zeta1}). The Riemann zeta function is strongly related to the prime numbers theory via the Euler product formula 
\begin{equation*}
\zeta (s) =\prod _{p} \left( 1-\frac{1}{p^ {s}}\right)^ {-1}
\end{equation*}
for $\Re s>1$, where $p$ ranges over primes.  The distribution of the zeta zeros is one of the outstanding problems in mathematics. We know that $\zeta (-2n) =0$ for every $n\geq 1$. The points $s=-2n$ are called the {\it trivial zeros} of the zeta function.  It is known that the possible {\it  non-trivial  zeros } of the zeta-function could only lie inside the {\it critical strip} $0<\Re s <1$. They are of great interest since their distribution leads to many important results in prime numbers theory.

We also know that the numbers of zeta  zeros is infinite and they lie symmetrical about the real axis and about the vertical  line $\Re s=\frac{1}{2}.$  {\it The Riemann hypothesis} posits  that all the non-trivial zeros lies on the {\it critical line } $\Re s =\frac{1}{2}. $

A probabilistic way to analyze the zeta zeros is to look at the values of $\log \zeta (s)$ on the critical line $s=\frac{1}{2} + \mathbf{i} t$ and to consider $t$ as a random variable uniformly distributed that takes huge values. That is, one considers $t\sim \mathcal{U}[T,2T]$ with $T$ close to infinity. By $\mathcal{U}[a,b]$ we will denote throughout this work  the uniform distribution over the interval $[a,b], a<b$.

Selberg' s theorem (see \cite{Se1}, \cite{Se2}, \cite{Se3} or the surveys \cite{K}, \cite{Tao}, \cite{Ts}) gives the asymptotic distribution of $\log \zeta (s) $ on the critical line $\Re s =\frac{1}{2}$.  Selberg theorem says that, if $t$ is a random variable uniformly distributed over the interval $[T, 2T]$, then the sequence
\begin{equation}
\label{logzeta}
\frac{ \log \zeta \left( \frac{1}{2}+ \mathbf{ i} t\right)}{\sqrt{\frac{1}{2}\log \log T}}
\end{equation}
converges in distribution to a complex-valued standard normal random variable $X_{1}+ \mathbf{i}X_{2}$
with $X_{1}, X_{2}\sim N(0,1)$ being independent random variables. There are several versions of this theorem. In particular, the result (\ref{logzeta}) holds if $t\sim U[0, T]$ of, more generally,  if $t\sim U[aT, bT]$ with $b>a \geq 0$.

Selberg's theorem basically says that the zeta zeros does not affect too much the behavior of  $\zeta$ on the critical line. Actually, the primes do most of the work.  The very small normalization of order $\sqrt{\log \log T }$ is usually interpreted as a {\it repulsion of zeros } (see \cite{Tao}, \cite{B}, \cite{CoDi}). 

Selberg theorem is actually  equivalent to the convergence of the real and imaginary parts of $\log \zeta $ on the critical line, i.e. (recall that, if $s$ is a complex number, then $\log s= \log \vert s \vert +\mathbf{i} \arg s$)
\begin{equation}
\label{s1}
\frac{  \log \left|\zeta \left( \frac{1}{2}+  \mathbf{i}t\right)\right|}{\sqrt{\frac{1}{2}\log \log T}}\to _{T\to \infty} ^{(d)} X_{1}\sim N(0,1)
\end{equation}
and
\begin{equation}
\label{s2}
\frac{ \arg \log \zeta \left( \frac{1}{2}+  \mathbf{i}t\right)}{\sqrt{\frac{1}{2}\log \log T}}\to _{T\to \infty} ^{(d)}X_{2} \sim N(0,1)
\end{equation}
where $" \to ^ {(d)}" $ stands for the convergence in distribution.  The idea of the proof of (\ref{s1}) and (\ref{s2}) is (see \cite{Se1}, \cite{Se2}, \cite{Se3}, \cite{Tao}, \cite{Titch}) to approximate $\log \zeta \left( \frac{1}{2}+\mathbf{i}t\right)$ by a Dirichlet series and  to look to the behavior of this Dirichlet series, which can be easier handled. The Dirichlet approximation will be of the form

\begin{equation}
\label{dir}
\sum _{p\leq T^ {\varepsilon} } \frac{1}{p ^{\frac{1}{2}+\mathbf{i}t}}=\sum _{p\leq T^ {\varepsilon} }\frac{\cos (t\log p)}{\sqrt{p}} + \mathbf{i}\sum _{p\leq T^ {\varepsilon} }\frac{\sin (t\log p)}{\sqrt{p}}.
\end{equation}
with $t\sim \mathcal{U}[aT, bT]$, $b>a \geq 0$ and   $\varepsilon$ small enough. We will work throughout with $\varepsilon =1$.

One of the main issues in our work is to study the speed of the convergence in (\ref{s1}) and (\ref{s2}). Some information concerning the rate of convergence to the normal distribution in the Selberg's central limit theorem can be found in the more recent Selberg's work \cite{Se4} or in \cite{Ts}, \cite{Wahl}. Actually, it follows from \cite{Se4} (see also Appendix A in \cite{Wahl} for a detailed proof) that (for $\varepsilon$ sutably chosen in (\ref{dir})) the Kolmogorov distance 
between the sequence (\ref{logzeta}) and the standard normal distribution is  less than $C\frac{\log \log \log T}{\sqrt{\log \log T}}$ (throughout, by $C$ we will denote a generic strictly positive constant that may change from one line to another).  It is actually shown in \cite{Se4}, \cite{Wahl} that the Kolmogorov distance between the Dirichlet series (\ref{dir})  and the  Gaussian law $N(0,1)$ is less $C\frac{1}{ \sqrt{\log \log T}}$  and it is deduced that  the  distance between (\ref{logzeta}) and $N(0,1)$ is less than  $C\frac{\log \log \log T}{\sqrt{\log \log T}}.$  It seems that there are not results concerning other metrics. Therefore, in a first step, we will  investigate the rate of convergence in the (one-dimensional)  Selberg's theorem in terms of the Wasserstein distance.

We use recent techniques based on Malliavin calculus combined with Stein method (see \cite{NPbook}) in order to obtain our error bounds. Let us describe our new contributions. First, concerning the one-dimensional Selberg's theorem: we prove that the distance (under several metrics, such as the Kolmogorov, total variation, Wasserstein or Fortet -Mourier metrics) between the Dirichlet approximation (\ref{dir}) and the standard normal law if  less than $C\frac{1}{\log \log T}$. This improves  the  known result  for the Kolmogorov distance. We also prove that the rate of convergence of $\log \zeta (\frac{1}{2}+ \mathbf{i} UT) $ to $N(0,1)$ is, under the Wasserstein metric, less than  $C\frac{1}{\sqrt{\log \log T}}$ and this improves the result in \cite{Se4}, \cite{Wahl} (which states, recall, that the Kolmogorov distance between (\ref{logzeta}) and $N(0,1)$ is less than  $C\frac{\log \log \log T}{\sqrt{\log \log T}}$). In our work, $U$ denotes  a standard uniform random variable  i.e. $U\sim \mathcal{U}[0,1]$.

We also study the multidimensional context. The multidimensional extension of the Selberg's theorem  has been proved more recently. First, in the paper \cite{HNY}, the authors showed that for any $0<\lambda _{1}<\ldots \lambda _{d}$, the random vector
\begin{equation}
\label{xt-intro} 
X_{T}:= \frac{1}{ \sqrt{\frac{1}{2} \log \log T}} \left( \log \zeta \left( \frac{1}{2} +\mathbf{i} P _{i}\right) \right) _{i=1,..,d}
\end{equation}
with $P_{i}= U e ^ {(\log T) ^ {\lambda _{i}}}$ ($i=1,\ldots , d$),  converges in distribution, as $T\to \infty$ to $(\lambda _{1} Y_{1}, \ldots , \lambda _{d} Y_{d})$ where $Y_{1}, \ldots , Y_{d}$ are independent standard complex Gaussian random variables.  There is no correlation between the components of the limit vector because the evaluation points $P_{i}$ are rather distant one from each other. When these points are less distanced, then non-trivial correlations appear in the limit. The result is due to \cite{B}. In this reference, the author showed that for $P_{i}= TU+ f_{T}^ {(i)}$ with $f_{T}^ {(i)}- f_{T} ^ {(j)}$ not too big (the exact meaning is given later), then the random vector (\ref{xt-intro}) converges in law to a $d$-dimensional complex Gaussian vector with dependent components. 

We  will regard these results from the Malliavin calculus point of view and we will give the associated error bounds.  We will  treat the case when $P_{i}= TU+ iT$ (here the space between points is pretty big and the limit is a Gaussian vector with independent components) and $P_{i} =TU+ f_{T}^ {(i)}$ with  $f_{T}^ {(i)}- f_{T} ^ {(j)}$ small if $i\not=j$ (non-trivial correlations appear). We will see that the order of the speed of convergence to the limit distribution is not affected by the distance between the evalution points. 

All these results lead to several consequences for the fluctuation of the number of zeros of the zeta function on the critical line.  More precisely, we prove that the number of zeta zeros on the critical line $\Re s =\frac{1}{2}$ between some random heights satisfies a central limit theorem and we obtain the associated error bound. Our results extend the findings in \cite{B}, \cite{CoDi} or \cite{HNY}.

Our paper is organized as follows.  In Section 2 we analyze the speed of convergence in the classical Selberg's theorem under several metrics via the Stein' s method combined with Malliavin calculus. In particular, we obtain explicit formulas for the Malliavin operators applied to the random variables in the left-hand side of (\ref{s1}), (\ref{s2}). In Section 3 we make the same study in the multidimensional settings, while in Section 4 we apply our findings to prove new results concerning the number of zeros of the Riemann zeta function. In the Appendix we included some elements from the Malliavin calculus and from the prime numbers theory needed in our work.

Throughout the paper we fix $H$ a real and separable Hilbert space and $(W(h), h\in H)$ an isonomal Gaussian process  (as introduced in Section \ref{mal})  on the probability space $\left(\Omega, \mathcal{F}, P\right)$.

\section{Rate of convergence  in the Selberg theorem via Malliavin calculus}

We here study the error bound corresponding to the weak convergences (\ref{s1}) and (\ref{s2}). The error bound will be obtained in two steps: first we measure the distance between the series  (\ref{dir}) and the standard normal distribution and then we will use an old result in \cite{Se3}.

\subsection{Rate of convergence for the Dirichlet series}

Consider the family $(X_{T}) _{T>0}$ given by
\begin{equation}
\label{xt}
X_{T}= \sum_{p\leq T}\left[ \frac{ \cos (TU\log p  )}{\sqrt{p}} - \mathbf{E}\frac{ \cos (TU\log p  )}{\sqrt{p}}\right]
\end{equation}
where the sum is taken over the primes $p$ and $U$ is $\mathcal{U}[0,1]$ distributed. In the sequel, we will assume 
\begin{equation}
\label{u}
U= e ^{-\frac{1}{2} (W(f)^{2}+W(g) ^{2}) } 
\end{equation}
 with $f,g \in H$, $\Vert f\Vert =\Vert g\Vert =1$ and $\langle f,g\rangle =0$ (all the scalar products and norms  in the paper will be considered in $H$ if no further precision is made).  In (\ref{u}), $W$ stands for a Gaussian isonormal process as described in the Appendix, Section \ref{mal}.  This implies that $W(f) $ and $W(g)$ are independent standard normal random variables.

The sequence $X_{T}$ corresponds to the real part of $\log \zeta \left( \frac{1}{2} + \mathbf{i} TU\right)$. A similar analysis can be done for the imaginary part. We briefly describe the main steps.

We will measure the distance between the sequence
$$\frac{1}{\sqrt{\frac{1}{2} \log \log T}} X_{T}$$
and the standard normal distribution. For every $p\leq T$, denote by
\begin{equation}
\label{xtp}
X _{T,p}= \frac{ \cos (TU\log p  )}{\sqrt{p}} - \mathbf{E}\frac{ \cos (TU\log p  )}{\sqrt{p}}
\end{equation}
so $X_{T}= \sum_{p\leq T  } X_{T,p}.$

We need to compute $ \langle DX_{T}, D(-L)^{-1} X _{T} \rangle$.  This quantity is crucial when one uses the Stein method combined with the Malliavin calculus (see \cite{NPbook}, see also (\ref{4i-4}) in Theorem \ref{t}). Let us do this computation.  We can write
\begin{eqnarray}
 &&\langle DX_{T}, D(-L)^{-1} X _{T} \rangle \nonumber\\
&=& \sum_{p_{1}, p_{2} \leq T} \langle DX_{T,p_{1}}, D(-L) ^{-1} X _{T, p_{2}}\rangle =\sum_{p_{1}, p_{2} \leq T} \langle DX_{T,p_{2}}, D(-L) ^{-1} X _{T, p_{1}}\rangle \nonumber\\
&=& \frac{1}{2} \sum_{p_{1}, p_{2} \leq T} \left( \langle DX_{T,p_{1}}, D(-L) ^{-1} X _{T, p_{2}}\rangle + \langle DX_{T,p_{2}}, D(-L) ^{-1} X _{T, p_{1}}\rangle \right).\label{21n-2}
\end{eqnarray}
Using the series expansion of the cosinus function $\cos x = \sum_{k\geq 0} \frac{(-1)^{k}} { (2k)! }x ^{2k}$ we get
$$X_{T,p}= \frac{1}{\sqrt{p} }\sum _{k\geq 0} \frac{(-1)^{k}} { (2k)! } (T\log p) ^{2k} (U ^{2k}-\mathbf{E} U ^{2k}).$$
So
\begin{eqnarray*}
&&\langle DX_{T,p_{1}}, D(-L) ^{-1} X _{T, p_{2}}\rangle +\langle DX_{T,p_{2}}, D(-L) ^{-1} X _{T, p_{1}}\rangle \\
&& = \sum_{k,l\geq 0} \frac{ (-1) ^{k+l}}{(2k)! (2l)!}T ^{2k+2l}
 \frac{ (\log p_{1}) ^{2k}}{\sqrt{p_{1}}}  \frac{ (\log p_{2}) ^{2l}}{\sqrt{p_{2}}}\\
&&\left[  \langle D (U^{2k}-\mathbf{E}U ^{2k}), D(-L) ^{-1} (U ^{2l}-\mathbf{E} U ^{2l})\rangle +\langle D (U^{2l}-\mathbf{E}U ^{2l}), D(-L) ^{-1} (U ^{2k}-\mathbf{E} U ^{2k})\rangle\right].
\end{eqnarray*}
Let us use  the notation, for every $k>0$
\begin{equation}
\label{gk}
G_{k}:= U ^{k}-\mathbf{E} U ^{k}=e ^{-\frac{k}{2}(W(f) ^{2}+ W(g) ^{2})}-\mathbf{E}e ^{-\frac{k}{2}(W(f) ^{2}+ W(g) ^{2})}.
\end{equation}
 Relation (\ref{21n-2}) becomes

\begin{eqnarray}
 \langle DX_{T}, D(-L)^{-1} X _{T} \rangle 
&=&\sum_{p_{1}, p_{2} \leq T}
\sum_{k,l\geq 0} \frac{ (-1) ^{k+l}}{(2k)! (2l)!}T ^{2k+2l}
\frac{ (\log p_{1}) ^{2k}}{\sqrt{p_{1}}} \frac{ (\log p_{2}) ^{2l}}{\sqrt{p_{2}}} \langle D G_k, D(-L) ^{-1}G_l\rangle \nonumber \\
&=& \frac{1}{2} \sum_{p_{1}, p_{2} \leq T} 
\sum_{k,l\geq 0} \frac{ (-1) ^{k+l}}{(2k)! (2l)!}T ^{2k+2l}
\frac{ (\log p_{1}) ^{2k}}{\sqrt{p_{1}}} \frac{ (\log p_{2}) ^{2l}}{\sqrt{p_{2}}} \nonumber \\
&& \left[
\langle D G_{2k}, D(-L) ^{-1} G_{2l} \rangle+\langle D G_{2l} , D(-L) ^{-1} G_{2k} \rangle\right].\nonumber\\
 \label{16n-4}
\end{eqnarray}
where we used  the symmetry of the sums. Consequently, it is necessary to calculate $\langle D G_{2k}, D(-L) ^{-1} G_{2l} \rangle+\langle D G_{2l} , D(-L) ^{-1} G_{2k} \rangle$. This will be done in the next lemma and it will  be used several times in the paper.

\begin{lemma}\label{l1}
Let $G_{k}$ be given by (\ref{gk}). Then for every $k,l> 0$
\begin{equation}
\langle D G_{2k}, D(-L) ^{-1} G_{2l} \rangle+\langle D G_{2l} , D(-L) ^{-1} G_{2k} \rangle =U ^{2k}  \frac{2k}{2l+1} \left( 1- U ^ {2l} \right)+ U ^{2l}  \frac{2l}{2k+1} \left( 1- U ^ {2k}\right)\label{15d-1}
\end{equation}
with $U$ given by (\ref{u}).

\end{lemma}
{\bf Proof: } 
For every two smooth centered random variables $F,G$ we have
\begin{eqnarray}
&&\langle DF, D(-L) ^{-1} G\rangle + \langle DG, D(-L) ^{-1} F\rangle \nonumber\\
&&=  \left( \langle D(F+G), D(-L) ^{-1} (F+G) \rangle - \langle DF, D(-L) ^{-1} F\rangle-\langle DG, D(-L) ^{-1} G\rangle \right)\label{pol}.
\end{eqnarray}
We use the following formula proved in \cite{NV}:
if  $Y= f(N)-{\mathbf E}[f(N)]$ where $f\in C_b^1( \mathbb{R} ^{n}; \mathbb{R})$ with bounded derivatives and $N=(N_{1},..., N_{n}) $ is a Gaussian vector   with zero mean and covariance matrix $K=(K_{i,j}) _{i,j=1,..,n}$ then
\begin{equation}\label{nv}
 \langle D(-L) ^{-1} (Y-{\mathbf E}[Y]), DY\rangle _H 
 = \int_0 ^1 da {\mathbf E}'\left[ \sum_{i,j=1} ^{n} K_{i,j} \frac{\partial f}{\partial x_{i}} (N) \frac{\partial f} {\partial x_{j}} (aN+ \sqrt{1-a^{2}}N' )\right].
\end{equation}
Here  $N'$ denotes an independent copy of $N$, the variables  $N$ and $N'$ are defined on a product probability space
 $\left( \Omega \times \Omega ', {\cal{F}} \otimes {\cal{F}}, P\times P'\right)$ and ${\mathbf E}'$
 denotes the expectation with respect to the probability measure $P'$.

In our case, for every $k\geq 1$,
$$U^{2k}-\mathbf{E}U^{2k}=G_{2k}= h( W(f), W(g))-\mathbf{E}h( W(f), W(g))$$
with $h(x,y)= e ^{-k(x ^{2}+ y^{2})}.$ Denote by 

\begin{equation}
\label{gka}
G_{k,a} = e ^{-k\left[ (aW(f) + \sqrt{1-a ^{2}}W'(f) ) ^{2}+ (aW(g) + \sqrt{1-a ^{2}}W'(g) ) ^{2}\right] }.
\end{equation}
Then, by (\ref{nv}) we find

\begin{eqnarray}
&&\langle D G_{2k}, D(-L) ^{-1} G_{2k}\rangle \nonumber \\
 &=& \int_{0} ^{1} da (2k) ^{2}  \mathbf{E}'\left[ U ^{2k}  G_{k,a}  W(f)( a(W(f)+ \sqrt{1-a^{2}} W'(f)) \right] \nonumber\\
&&+ \int_{0} ^{1} da (2k) ^{2}  \mathbf{E}'\left[ U ^{2k}  G_{k,a}  W(g)( a(W(g)+ \sqrt{1-a^{2}} W'(g)) \right] \nonumber \\
&=&\int_{0} ^{1} da 4k^{2}a  U ^{2k} \left( W(f) ^{2} + W(g) ^{2}\right) \mathbf{E}'(G_{k,a}) \nonumber \\
&&+ \int_{0} ^{1} da 4k^{2} \sqrt{1- a ^{2}}U ^{2k} \left[ W(f) \mathbf{E}' (W'(f) G_{k,a}) + W(g) \mathbf{E}'(W'(g) G_{k, a} ) \right].\label{16n-2}
\end{eqnarray}
Let us first calculate $\mathbf{E}'(G_{k,a})$  with $G_{k,a}$ given by (\ref{gka}).  We have
\begin{equation}
\label{16n-1}\mathbf{E}'(G_{k, a} )= g(aW(f)) g(aW(g))
\end{equation}
where $g(c)= \mathbf{E}  e ^{ -k( c+ \sqrt{1- a ^{2}} Z ) ^{2}}$ with $Z$ a standard normal random variable. By standard calculations
\begin{eqnarray*}
g(c)&=& \frac{1}{\sqrt{2\pi}} \int_{\mathbb{R}} e ^{ -k( c+ \sqrt{1- a ^{2}} x ) ^{2} }e ^{-\frac{x^{2}}{2}}dx 
=\frac{1}{\sqrt{2\pi}} e ^{-k c^{2}}\int_{\mathbb{R}}dx e ^{-k(1-a ^{2} ) x ^{2} }e ^{-\frac{x^{2}}{2}} e ^{-2kc \sqrt{1- a ^{2}}x }\\
&=&\frac{1}{\sqrt{2\pi}} e ^{-k c^{2}}  e ^{ \frac{ 2k ^{2} c ^{2} (1-a ^{2})}{1+2k(1-a ^{2})}}\int_{\mathbb{R}} e ^{-\frac{1}{2} (1+2k(1-a ^{2}) ) y^{2} }dy =\frac{1}{ \sqrt{1+2k(1-a ^{2})} } e ^{-\frac{kc^{2}}{1+2k(1-a ^{2})}}.
\end{eqnarray*}
By (\ref{16n-1}), we obtain
\begin{equation}
\label{aux1}\mathbf{E}' (G_{k,a})= \frac{1}{ 1+2k(1-a^{2}) } e ^{-\frac{k a^ {2}(W(f) ^{2} + W(g) ^{2}) }{ 1+2k(1-a ^{2}) }}.
\end{equation}
We also need to compute $\mathbf{E}'(G_{k,a} W '(f).$ We have
$$\mathbf{E}'(G_{k,a} W '(f))= \mathbf{E}'\left[ W'(f) e ^{-k(aW(f)+ \sqrt{1- a ^{2}}W'(f) ) ^{2} }\right] g(a W(g))$$
where $g(c)$ has been computed just above. We will find that 

$$\mathbf{E}'\left[ W'(f) e ^{-k(aW(f)+ \sqrt{1- a ^{2}}W'(f) ) ^{2} }\right]  = m(aW(f))$$
with $m(c)= \mathbf{E}' \left[ Z  e ^{-k (c+ \sqrt{1-a ^{2}} Z ) ^{2} }\right].$ Moreover
\begin{eqnarray*}
m(c)&=& \frac{1}{\sqrt{2\pi}} \int_{\mathbb{R}} xe ^{ -k( c+ \sqrt{1- a ^{2}} x ) ^{2} }e ^{-\frac{x^{2}}{2}}dx \\
&=& e ^{-\frac{kc^{2}}{1+2k(1-a ^{2})}}\frac{1}{\sqrt{2\pi}} 
\int_{\mathbb{R}} (y- \frac{2kc\sqrt{1-a ^{2}}}{ 1+2k (1-a ^{2})  }) e ^{-\frac{1}{2} (1+2k(1-a ^{2}) ) y^{2} }dy \\
&=&-e ^{-\frac{kc^{2}}{1+2k(1-a ^{2})}} \frac{2kc\sqrt{1-a ^{2}}}{ (1+2k (1-a ^{2})) ^{\frac{3}{2}} }.
\end{eqnarray*}
Hence
\begin{equation}
\label{aux2}
\mathbf{E}'(G_{k,a} W '(f)) =-e ^{-\frac{k a^{2}  (W(f)^{2} +W(g) ^{2})} {1+2k(1-a ^{2})}} \frac{2ka W(f)\sqrt{1-a ^{2}}}{ (1+2k (1-a ^{2})) ^{2} }
\end{equation}
and similarly,

\begin{equation}
\label{aux3}
\mathbf{E}'(G_{k,a} W '(g)) =-e ^{-\frac{k a^{2} (W(f)^{2} +W(g) ^{2})} {1+2k(1-a ^{2})}} \frac{2ka W(g)\sqrt{1-a ^{2}}}{ (1+2k (1-a ^{2})) ^{2} }.
\end{equation} 
Now, formula (\ref{16n-2}) becomes

\begin{eqnarray*}
&&\langle D G_{2k}, D(-L) ^{-1} G_{2k}\rangle\\
&=& \int_{0} ^{1} da 4 k^{2} a   \frac{1}{1+2k(1-a^{2})}U^{2k} ( W(f) ^{2} + W(g) ^{2}) e ^{-\frac{k a^{2} (W(f)^{2} +W(g) ^{2})} {1+2k(1-a ^{2})}} \\
&&-\int_{0} ^{1} da 4 k^{2} (1-a ^{2 } ) \frac{2ka }{ (1+2k (1-a ^{2})) ^{2} } U^{2k} ( W(f) ^{2} + W(g) ^{2}) e ^{-\frac{k a^{2} (W(f)^{2} +W(g) ^{2})} {1+2k(1-a ^{2})}}.
\end{eqnarray*}
Let us denote by
\begin{equation}
\label{s}
S=W(f)^{2} +W(g) ^{2}.
\end{equation}
Since
\begin{equation}
\label{aux4}a \frac{1}{1+2k(1-a^{2})}- (1- a^{2}) \frac{2ka }{ (1+2k (1-a ^{2})) ^{2} }= \frac{a}{ (1+2k (1-a ^{2})) ^{2} }
\end{equation}
we can write
\begin{eqnarray*}
\langle D G_{2k}, D(-L) ^{-1} G_{2k}\rangle &=&4k^{2}U ^{2k} S \int_{0} ^{1} da \frac{a}{ (1+2k (1-a ^{2})) ^{2} }e ^{-\frac{k a^{2} S} {1+2k(1-a ^{2})}}.
\end{eqnarray*}
Note that

\begin{equation}
\label{aux5} \left( \frac{a^{2}}{ (1+2k (1-a ^{2})) ^{2} }\right) '= \frac{2a(1+2k)}{ (1+2k (1-a ^{2})) ^{2} }.
\end{equation}
Thus, by the change of variables  $\frac{a^{2}}{ (1+2k (1-a ^{2})) ^{2} }=z$
\begin{eqnarray*}
\langle D G_{2k}, D(-L) ^{-1} G_{2k}\rangle &=&4k^{2}U ^{2k} S \frac{1}{2(2k+1)} \int_{0} ^{1} dz e ^{-kSz} \\
&=&4k^{2}U ^{2k} S \frac{1}{2(2k+1)} \int_{0} ^{1} dz e ^{-kSz} =2kU ^{2k}  \frac{1}{2k+1} \left( 1- e ^{-kS} \right).
\end{eqnarray*}
We  can compute
$$ \langle D(G_{2k}+ G_{2l}), D(-L) ^{-1} (G_{2k}+ G_{2l}) \rangle$$
by using again (\ref{nv}) with $h(x,y)= e ^{-k(x ^{2}+ y^{2})}+ e ^{-l(x ^{2}+ y^{2})}.$
We will have 
\begin{eqnarray*}
&&\langle D(G_{2k}+ G_{2l}), D(-L) ^{-1} (G_{2k}+ G_{2l}) \rangle \\
&=& \int_{0} ^{1} da a (2kU^{2k}+ 2lU^{2l}) (W(f) ^{2} + W(g) ^{2}) \mathbf{E}'( 2kG_{k,a} + 2l G_{l,a} ) \\
&&+ \int_{0} ^{1} da \sqrt{1- a ^{2}} (2kU^{2k}+ 2lU^{2l})  \\
&&\left[ W(f) \mathbf{E}'W'(f) (2kG_{k,a}+ 2lG_{l,a} )  + W(g) \mathbf{E}'W'(g) (2kG_{k,a}+ 2lG_{l,a}) \right]
\end{eqnarray*}
and this implies
\begin{eqnarray*}
&&\langle D(G_{2k}+ G_{2l}), D(-L) ^{-1} (G_{2k}+ G_{2l}) \rangle - \langle DG_{2k}, D(-L) ^{-1} G_{2k} \rangle - \langle DG_{2l}, D(-L) ^{-1} G_{2l} \rangle\\
&=& \int_{0} ^{1} da a 4kl (W(f) ^{2} + W(g) ^{2}) \left[ U^{2k} \mathbf{E}'G_{l,a} +U^{2l} \mathbf{E} '(G_{k,a} \right] \\
&&+\int_{0} ^{1} da \sqrt{1- a ^{2}}  4kl \left[ U^{2k} W(f) \mathbf{E}'W'(f) G_{l,a}+U^{2l} W(f) \mathbf{E}'W'(f) G_{k,a}\right.\\
&&\left.  +U^{2k} W(g) \mathbf{E} 'W'(g) G_{l,a} + U^ {2l} W(g) \mathbf{E} 'W'(g) G_{k,a} \right].
\end{eqnarray*}
Consequently, from relations (\ref{aux1}), (\ref{aux2}), and (\ref{aux3}), we obtain
\begin{eqnarray*}
&&\langle D(G_{k}+ G_{l}), D(-L) ^{-1} (G_{k}+ G_{l}) \rangle - \langle DG_{k}, D(-L) ^{-1} G_{k} \rangle - \langle DG_{l}, D(-L) ^{-1} G_{l} \rangle\\
&=& \int_{0} ^{1} da a 4kl (W(f) ^{2} + W(g) ^{2}) \\
&&\left[ \frac{1}{1+k(1-a^{2})} G_{l}e ^{-\frac{k a^{2} (W(f)^{2} +W(g) ^{2})} {1+2k(1-a ^{2})}} + \frac{1}{1+2l(1-a^{2})} G_{k}e ^{-\frac{l a^{2} (W(f)^{2} +W(g) ^{2})} {1+2l(1-a ^{2})}}\right] \\
&&-\int_{0} ^{1} da \sqrt{1- a ^{2}}  4kl  (W(f)^{2} +W(g) ^{2})\\
&&\left[ G_{l} e ^{-\frac{k a^{2}  (W(f)^{2} +W(g) ^{2})} {1+2k(1-a ^{2})}} \frac{2ka \sqrt{1-a ^{2}}}{ (1+2k (1-a ^{2})) ^{2} }  + G_{k} e ^{-\frac{l a ^{2}(W(f)^{2} +W(g) ^{2})} {1+2l(1-a ^{2})}} \frac{2la \sqrt{1-a ^{2}}}{ (1+2l (1-a ^{2})) ^{2} }\right] .
\end{eqnarray*}
To conclude (\ref{15d-1}), it suffices to use (\ref{pol}),   (\ref{aux4}) and (\ref{aux5}). \qed

We  obtain the explicit form of the terms needed in the Stein-Malliavin bound (\ref{4i-4}).

\begin{prop}\label{p1}
For every $T>0$, with $X_{T}$ given by (\ref{xt}), we have
\begin{eqnarray}\label{15d-5}
&&\langle DX_{T}, D(-L)^{-1} X _{T} \rangle =  \sum_{p_{1}, p_{2} \leq T}\frac{1}{\sqrt{p_{1}p_{2}}}\frac{\log p_{2} }{\log p_{1}} \nonumber \\
& &  \times \left[\sin (TU\log p_{2} )\sin (TU\log p_{1} ) -
 U\sin (TU\log p_{2} )\sin (T \log p_{1}) \right].
\end{eqnarray}

\end{prop}
{\bf Proof: } By (\ref{16n-4}) and relation (\ref{15d-1}) in Lemma \ref{l1}, we have

\begin{eqnarray}
 \langle DX_{T}, D(-L)^{-1} X _{T} \rangle  
&=&\frac12\sum_{p_{1}, p_{2} \leq T}
\sum_{k,l\geq 0} \frac{ (-1) ^{k+l}}{(2k)! (2l)!}T ^{2k+2l}
\frac{ (\log p_{1}) ^{2k}}{\sqrt{p_{1}}} \frac{ (\log p_{2}) ^{2l}}{\sqrt{p_{2}}} \nonumber \\
&& \times \left[2kU ^{2k}  \frac{1}{2l+1} \left( 1- e ^{-lS} \right)+ 2lU ^{2l}  \frac{1}{2k+1} \left( 1- e ^{-kS} \right)\right] \nonumber \\
&=&\sum_{p_{1}, p_{2} \leq T}
\sum_{k\geq 0; l\geq 1} \frac{ (-1) ^{k+l}}{(2k+1)! (2l-1)!}T ^{2k+2l}
\frac{ (\log p_{1}) ^{2k}}{\sqrt{p_{1}}} \frac{ (\log p_{2}) ^{2l}}{\sqrt{p_{2}}} U ^{2l}  \left( 1- e ^{-kS} \right).\nonumber
\end{eqnarray}
So,
\begin{eqnarray}
 \langle DX_{T}, D(-L)^{-1} X _{T} \rangle 
&=&\sum_{p_{1}, p_{2} \leq 
T}\frac{1}{\sqrt{p_{1}p_{2}}}\sum_{l \geq 1} \frac{(-1) ^{l}}{(2l-1) !} (T \log p_{2}) ^{2l} U^{2l}\nonumber\\
&& \sum_{k\geq 0} \frac{(-1) ^{k} }{(2k+1) !} (T\log p_{1}) ^{2k} (1- U ^{2k}).\label{15d-2}
\end{eqnarray}
We first compute the sum over $l$. We have
\begin{eqnarray}
\sum_{l \geq 1} \frac{(-1) ^{l}}{(2l-1) !} (T \log p_{2}) ^{2l} G_{l}&=& \sum_{l\geq 0} \frac{(-1) ^{l+1}}{(2l+1) !} (T \log p_{2}) ^{2l+2} U ^{2l+2}\nonumber\\
&=& -TU\log p_{2} \sin (TU\log p_{2} ). \label{15d-3}
\end{eqnarray}
Concerning the sum over $k$
\begin{eqnarray}
&&\sum_{k\geq 0} \frac{(-1) ^{k} }{(2k+1) !} (T\log p_{1}) ^{2k} (1- U ^{2k})
=\sum_{k\geq 0} \frac{(-1) ^{k} }{(2k+1) !} (T\log p_{1}) ^{2k} (1- U ^{2k})\nonumber\\
&=&\frac{1}{T\log p_{1}} 
\left[ \sum_{k\geq 0} \frac{(-1) ^{k} }{(2k+1) !} (T\log p_{1}) ^{2k+1}  - U ^{-1}\sum_{k\geq 0} \frac{(-1) ^{k} }{(2k+1) !} (T\log p_{1}) ^{2k+1} U ^{2k+1}\right] \nonumber \\
&=&\frac{1}{T\log p_{1}}  \left( \sin (T \log p_{1}) - U ^{-1}\sin (TU\log p_{1} ) \right).\label{15d-4}
\end{eqnarray}
By combining (\ref{15d-2}), (\ref{15d-3}) and (\ref{15d-4}), we get (\ref{15d-5}). \qed

\vskip0.2cm

Using the same lines, we can treat the imaginary part of $\log \zeta (\frac{1}{2}+\mathbf{i}t)$, with $t=TU$.

\begin{prop}
Denote by 
\begin{equation}
\label{yt}
Y_{T} = \sum_{p\leq T}\left( \frac{\sin (TU \log p  )}{\sqrt{p}}-\frac{\mathbf{E}\sin (TU \log p  )}{\sqrt{p}}\right). 
\end{equation}
Then
\begin{eqnarray*}
\langle DY_{T}, D(-L) ^{-1} Y_{T}\rangle &=& \sum_{p_{1}, p_{2} \leq T}\frac{\log p_{1}}{\log p_{2}} \frac{1}{\sqrt{p_{1}p_{2}}}\left[ (U-1) \cos (TU\log p_{2})-U \cos (TU\log p_{2}) \cos (T\log p_{1}) \right.\\
&&\left. + \cos (TU\log p_{1} ) \cos (TU\log p_{2} )\right]. 
\end{eqnarray*}
\end{prop}
{\bf Proof: } As in the proof of Proposition \ref{p1}, 
\begin{eqnarray}
 &&\langle DY_{T}, D(-L)^{-1} Y _{T} \rangle \nonumber \\
&=& \frac{1}{2} \sum_{p_{1}, p_{2} \leq T} \frac{1}{\sqrt{p_{1}p_{2}}}
\sum_{k,l\geq 0} \frac{ (-1) ^{k+1++1l}}{(2k+1)! (2l+1)!}T ^{2k+1+2l+1}
\frac{ (\log p_{1}) ^{2k+1}}{\sqrt{p_{1}}} \frac{ (\log p_{2}) ^{2l+1}}{\sqrt{p_{2}}} \nonumber \\
&&\left[ \langle D G_{2k+1}, D(-L) ^{-1} G_{2l+1} \rangle-\langle D G_{2l+1} , D(-L) ^{-1} G_{2k+1} \rangle\right] \nonumber\\
\end{eqnarray}
with $G_{2k+1}$ from (\ref{gk}) and by Lemma \ref{l1} (by replacing $k,l$ by $k+\frac{1}{2}, l+\frac{1}{2}$ respectively),
\begin{eqnarray*}
&& \left[ \langle D G_{2k+1}, D(-L) ^{-1} G_{2l+1} \rangle-\langle D G_{2l+1} , D(-L) ^{-1} G_{2k+1} \rangle\right] \\
&=& (2k+1)U ^{2k+1}  \frac{1}{2l+2} \left( 1- U ^{2l+1} \right)+ (2l+1)U ^{2l+1}  \frac{1}{2k+2} \left( 1- U ^{2k+1} \right).
\end{eqnarray*}
Thus
\begin{eqnarray*}
&&\langle DY_{T}, D(-L)^{-1} Y _{T} \rangle  \\
&=& \sum_{p_{1}, p_{2} \leq T}  \frac{1}{\sqrt{p_{1}p_{2}}}\sum_{k\geq 0} \frac{(-1) ^{k}}{(2k+2)!} (T\log p_{1} ) ^{2k+1} (1- U ^{2k+1}) \sum _{l\geq 0} \frac{(-1) ^{l}}{(2l)!} (TU\log p_{2}  ) ^{2l+1}.
\end{eqnarray*}
To conclude, it remains to notice that
$$\sum _{l\geq 0} \frac{(-1) ^{l}}{(2l)!} (TU\log p_{2}  ) ^{2l+1}=TU \log p_{2} \cos (TU \log p_{2} )$$
and
$$\sum_{k\geq 0} \frac{(-1) ^{k}}{(2k+2)!} (T\log p_{1} ) ^{2k+1} (1- U ^{2k+1})= \frac{1}{T\log p_{1}}\left( 1-\cos (T\log p_{1})+ \frac{1}{U} (\cos (T\log p_{1}U)-1)\right).$$
\qed

Let us take a moment to introduce certain notion on the distance between  probability distributions and to recall some links between these topic and Malliavin calculus. Let $X, Y$ be two  random variables. 
The distance between the law of  $X$ and the law of $Y$ is usually defined by (${\mathcal{L}}(F)$  denotes the law of $F$)
 \begin{equation*}
 d({\cal{L}}(X), {\cal{L}}(Y)): = \sup _{h\in {\cal{H}} } \vert \mathbf{E} h(X)- \mathbf{E}h(Y) \vert
 \end{equation*}
where ${\cal{H}} $ is a suitable class of functions. For example, if ${\cal{H}} $  is the set of indicator functions
 $1_{(-\infty, z]} , z\in \mathbb{R}$  we obtain  the Kolmogorov \index{Kolmogorov} distance (for simplicity, we will always write $d(X,Y)$ instead of   $d({\cal{L}}(X), {\cal{L}}(Y))$)
 \begin{equation*}
 d_{K}(X,Y):=d_{K} ({\cal{L}}(X),{\cal{L}}(Y))= \sup_{z \in \mathbb {R}} \vert P (X\leq z) -P (Y\leq z) \vert .
 \end{equation*}
 If ${\cal{H}} $ is  the set of $1_{B}$ with $B$ a Borel set, one has the total variation  distance
 \begin{equation*}
d_{TV} ({\cal{L}}(X),{\cal{L}}(Y)) =\sup_{B\in {\cal{B}}(\mathbb{R}) }\left|P(X\in B)-P(Y\in B) \right|
\end{equation*}
while for  ${\cal{H}}=\{ h; \Vert h\Vert _{L} \leq 1\} $ ($\Vert \cdot \Vert _{L}$ is the Lipschitz norm) one has the Wasserstein  distance denoted $d_{W}$. We will focus in our work on these metrics.  We will use the generic notation $d(X,Y)$ when our claim concerns all the metrics introduced above.

Let us recall the Stein bound for the normal approximation in terms of the Malliavin operators. See Section 5 in \cite{NPbook}.

\begin{theorem}\label{t}
If $F$ is a random variable in $\mathbb{D} ^{1,4}$ with $\mathbf{E}F=0$ and $N$ is a standard normal random variable, then
\begin{equation}
\label{4i-4}d (F, N)\leq C \mathbf{E}\left| 1-\langle DF, D(-L) ^{-1} F \rangle \right|. 
\end{equation}
\end{theorem}

We have the following result.

\begin{theorem}\label{t1}
For every $T>0$, let $X_{T}, Y_{T}$ be given by (\ref{xt}), (\ref{yt}) respectively. Denote by
\begin{equation}
\label{ft}
F_{T}= \frac{1}{\sqrt {\frac{1}{2}\log \log T}}X_{T}, \hskip0.3cm G_{T}=\frac{1}{\sqrt {\frac{1}{2}\log \log T}}Y_{T}.
\end{equation}
Then for $T$ large enough,
\begin{equation}
\label{4i-1}d (F_{T}, N) \leq C\frac{1}{ \log \log T } \mbox{ and } d (G_{T}, N) \leq C\frac{1}{ \log \log T } .
\end{equation}
\end{theorem}
{\bf Proof: }Clearly  $\mathbf{E}F_{T}=0$ and
from relation (\ref{15d-5}) in Proposition \ref{p1}

\begin{eqnarray*}
 &&\langle DX_{T}, D(-L)^{-1} X _{T} \rangle  \\
&=&- \sum_{p_{1}, p_{2} \leq T}\frac{1}{\sqrt{p_{1}p_{2}}}\frac{\log p_{2} }{\log p_{1}}U\sin (TU\log p_{2} )
\left( \sin (T \log p_{1}) - U ^{-1 }\sin (TU\log p_{1} ) \right)\\
&=&-  \sum_{p_{1}, p_{2} \leq T}\frac{1}{\sqrt{p_{1}p_{2}}}\frac{\log p_{2} }{\log p_{1}}U\sin (TU\log p_{2} )
\sin (T \log p_{1}) \\
&&+ \sum_{p_{1}, p_{2} \leq T}\frac{1}{\sqrt{p_{1}p_{2}}}\frac{\log p_{2} }{\log p_{1}} \sin (TU\log p_{2} ) \sin (TU\log p_{1} )
\end{eqnarray*}
and by separating the diagonal and non-diagonal parts in the second sum above, and by using $\sin ^ {2}(x)= \frac{1-\cos(2x)}{2}$,
\begin{eqnarray*}
 &&\langle DX_{T}, D(-L)^{-1} X _{T} \rangle  \\
&=& \sum _{p\leq T} \frac{\sin^{2} (TU \log p )}{p} + \sum_{p_{1}, p_{2} \leq T; p_{1}\not=p_{2} }\frac{1}{\sqrt{p_{1}p_{2}}}\frac{\log p_{2} }{\log p_{1}} \sin (TU\log p_{2} ) \sin (TU\log p_{1} )\\
&&-  \sum_{p_{1}, p_{2} \leq T}\frac{1}{\sqrt{p_{1}p_{2}}}\frac{\log p_{2} }{\log p_{1}}U\sin (TU\log p_{2} )
\sin (T \log p_{1})\\
&=&\sum _{p\leq T} \frac{1}{2p} -\sum _{p\leq T} \frac{\cos (2T U\log p )}{2p}
+ \sum_{p_{1}, p_{2} \leq T; p_{1}\not=p_{2} }\frac{1}{\sqrt{p_{1}p_{2}}}\frac{\log p_{2} }{\log p_{1}} \sin (TU\log p_{2} ) \sin (TU\log p_{1} )\\
&&-  \sum_{p_{1}, p_{2} \leq T}\frac{1}{\sqrt{p_{1}p_{2}}}\frac{\log p_{2} }{\log p_{1}}U\sin (TU\log p_{2} )
\sin (T \log p_{1}).
\end{eqnarray*}
Therefore, 

\begin{eqnarray*}
\mathbf{E} \left| 1-\langle DF_{T}, D(-L) ^{-1} F_{T} \rangle \right| &\leq & A_{1,T}+ A_{2,T}+ A_{3,T}+ A_{4,T}+ A_{5,T}
\end{eqnarray*}
where

\begin{equation*}
A_{1,T}=\left| 1- \frac{1}{\frac{1}{2}\log \log T}\sum _{p\leq T} \frac{1}{2p} \right|, A_{2,T}= \frac{1}{\frac{1}{2}\log \log T} \left|  \mathbf{E}\sum _{p\leq T} \frac{2\cos(2T U\log p )}{2p}\right|,
\end{equation*}

$$A_{3,T}= \frac{1}{\frac{1}{2}\log \log T}\left| \mathbf{E}\sum_{p_{1}, p_{2} \leq T ^{\varepsilon}; p_{1}\not=p_{2} }\frac{1}{\sqrt{p_{1}p_{2}}}\frac{\log p_{2} }{\log p_{1}} \sin (TU\log p_{2} ) \sin (TU\log p_{1} )\right| $$
and
$$A_{4,T}= \frac{1}{\frac{1}{2}\log \log T}\left| \mathbf{E}\sum_{p_{1}, p_{2} \leq T}\frac{1}{\sqrt{p_{1}p_{2}}}\frac{\log p_{2} }{\log p_{1}}U\sin (TU\log p_{2} )
\sin (T \log p_{1})\right|.$$
Since for every $a \in \mathbb{R}, $ one has $ \mathbf{E} \cos ( a T U\log p ) =\frac{1}{a T\log p} \sin (a T \log p)$, we have 
\begin{equation}
\label{15d-9}\left|   \mathbf{E} \cos ( a TU \log p )  \right| \leq C \frac{1}{T \log p}.
\end{equation}
The inequality (\ref{15d-9}) gives immediately, via (\ref{pn2})
\begin{equation*}
\vert A_{2,T} \vert \leq C \frac{1}{T \log \log T}\sum _{p\leq T} \frac{1}{p \log p}\leq C \frac{1}{T \log \log T}\sum _{p\leq T} \frac{1}{p }\leq C\frac{1}{T}.
\end{equation*}
To bound $A_{3,T}$,  we use $\sin  (x) \sin (y) =\frac{1}{2} (\cos (x-y)- \cos (x-y)) $ and, if $p_{1}\not=p_{2}$,  
\begin{equation}
\label{15d-6}
\left| \mathbf{E} \cos (aTU(\log p_{1}\pm \log p_{2})\right| \leq C \frac{1}{T}
\end{equation}
to get
\begin{eqnarray*}
\vert A_{3,T} \vert &\leq& C \frac{1}{T \log \log T} \sum _{p_{1}\leq T} \frac{1}{\sqrt{p_1}\log p_{1}}\sum _{p_{2}\leq 
T} \frac{\log p_{2}}{\sqrt{p_2}}\leq C\frac{1}{T \log \log T} \sum _{p_{1}\leq T} \frac{1}{\sqrt{p_1}}\log T \sum _{p_{2}\leq T} \frac{1}{\sqrt{p_2}}\\
&\leq& C
\frac{1}{ \log \log T \log T}\frac{1}{ \log  \log T}
\end{eqnarray*}
where we used the estimate (\ref{pn1}). Finally, to deal  with the summand $A_{4,T}$, we majorize $\vert \sin (T \log p_{1})\vert $ by 1, we use
 \begin{equation}
\label{18d-2}
\vert \mathbf{E}U\sin (TU\log p_{2} )\vert \leq C \frac{1}{T\log p_{2}}
\end{equation} and we will have, from (\ref{18d-2}) and (\ref{pn1})
\begin{equation*}
\vert A_{4,T} \vert \leq C \frac{1}{T \log \log T} \sum _{p_{1}\leq T} \frac{1}{\sqrt{p_1}\log p_{1}}\sum _{p_{2}\leq T} \frac{1}{\sqrt{p_2}}\leq C \frac{1}{ \log \log T (\log T) ^{2}}.
\end{equation*}
The speed of the convergence will be given by the dominant term $A_{1,T}$. Actually, from (\ref{pn2}), 
$$\vert A_{1,T} \vert \leq C \frac{1}{\log \log T}.$$
Obviously, similar arguments apply to the sequence $G_{T}$ from (\ref{ft}). \qed

\begin{remark}
Our inequalities (\ref{4i-1}) improve the bounds obtained in \cite{Se4} (see also Appendix A in \cite{Wahl}) where it was proved that for large $T$
$$d_{K} (F_{T}, N) \leq C\frac{1}{\sqrt{ \log \log T }} \mbox{ and } d _{K}(G_{T}, N) \leq C\frac{1}{ \sqrt{\log \log T} } $$
where $d_{K}$ is the Kolmogorov distance.
\end{remark}

\subsection{Rate of convergence in the Selberg theorem}

The real part of $\log \zeta (s)$ on the critical, where $\zeta$ is the Riemann zeta function (\ref{zeta1}) , can be approximated by the family $X_T + \mathbf{E}X_T$ where $X_T$ is given by (\ref{xt}). More precisely, if $t$ is a random variable uniformly distributed on $[0,T]$, then $\log \left| \zeta (\frac{1}{2} + \mathbf{i}t ) \right| $ is "close" (we explain below what that means) to $\sum _{p\leq T} \frac{ \cos (TU\log p ) }{\sqrt{p}}$ with $U\sim \mathcal{U}[0,1]$.  Since we have estimated is the previous paragraph the distance between $X_{T}$ and the standard normal law, we will be able to measure how far is $\log \left| \zeta (\frac{1}{2} + \mathbf{i}t ) \right| $ from the standard normal distribution. Actually, we have (with $t=TU\sim \mathcal{U}[0,T]$)

\begin{eqnarray*}
\frac{\log \vert \zeta \left( \frac{1}{2}+\mathbf{i}t\right) \vert}{\sqrt{\frac{1}{2} \log \log T}} 
&=& 
\frac{1}{ \sqrt{\frac{1}{2}\log \log T}}\sum_{p\leq T} \frac{ \cos (TU\log p  )- \mathbf{E} \cos (TU\log p ) }{\sqrt{p}}\\
&+& \frac{1}{ \sqrt{\frac{1}{2}\log \log T}}\left( \log \vert \zeta \left( \frac{1}{2}+\mathbf{i}t\right) \vert -\sum_{p\leq T} \frac{ \cos (TU\log p  ) }{\sqrt{p}}\right) \\
&+& \frac{1}{\sqrt{ \frac{1}{2}\log \log T}}\sum_{p\leq T} \frac{\mathbf{E} \cos (TU\log p  ) }{\sqrt{p}}.
\end{eqnarray*}

As mentioned above, the second summand in the right side above is known to be "small". The exact meaning is described in the  following result which has been obtained in \cite{Se3}. The reader may also consult the survey \cite{K} for the detailed steps of the proof.

\begin{lemma} \label{l2}For every $k,j\geq 1$ integers,  it holds 
\begin{equation}\label{15d-7}
\mathbf{E} \left| \log \vert \zeta \left( \frac{1}{2}+\mathbf{i}T(U+j)\right) \vert -\sum_{p\leq T} \frac{ \cos (T\log p (U+j) ) }{\sqrt{p}}\right| ^{2k}= \mathcal{O} (1).
\end{equation}

\end{lemma}

We do not need to assume the Riemann hypothesis in order to have the result in Lemma \ref{l2}. 
We will use the Wasserstein distance (introduced in this section) to measure how far is $\log \left| \zeta (\frac{1}{2} + \mathbf{i}t ) \right| $ from the standard normal distribution.  From the definition of this metric, one can see that
\begin{equation}
\label{15d-8}d _W (F,G) \leq \mathbf{E} \vert F-G \vert \leq \left( \mathbf{E}\left| F-G\right| ^{2} \right) ^{\frac{1}{2}}
\end{equation}
if $F,G$ are two random variables in $L^2 (\Omega)$.  Using the triangular inequality for the Wasserstein distance, we write, with  $t=TU$ and $U$ as in (\ref{u})

\begin{eqnarray}
&& d_{W} \left( \frac{\log \vert \zeta \left( \frac{1}{2}+\mathbf{i}t\right) \vert}{\sqrt{\frac{1}{2} \log \log T}} , N(0,1)\right)\nonumber\\
&\leq &d_{W} \left( \frac{1}{ \sqrt{\frac{1}{2}\log \log T}}\sum_{p\leq T} \frac{ \cos (TU\log p  )- \mathbf{E} \cos (TU\log p ) }{\sqrt{p}}, N(0,1)\right)\nonumber \\
&+&d_{W} \left( \frac{\log \vert \zeta \left( \frac{1}{2}+\mathbf{i}t\right) \vert}{\sqrt{\frac{1}{2} \log \log T}}, \frac{1}{ \sqrt{\frac{1}{2}\log \log T}}\sum_{p\leq T} \frac{ \cos (TU\log p  )- \mathbf{E} \cos (TU\log p ) }{\sqrt{p}}\right)\nonumber\\
&\leq & d_{W} \left( \frac{1}{ \sqrt{\frac{1}{2}\log \log T}}\sum_{p\leq T} \frac{ \cos (TU\log p  )- \mathbf{E} \cos (TU\log p ) }{\sqrt{p}}\right)\nonumber \\
&+& \frac{1}{\sqrt{ \frac{1}{2}\log \log T}}\mathbf{E} \left| \log \vert \zeta \left( \frac{1}{2}+\mathbf{i}t\right) \vert -\sum_{p\leq T} \frac{ \cos (TU\log p  ) }{\sqrt{p}}\right|\nonumber \\
&+& \frac{1}{\sqrt{ \frac{1}{2}\log \log T}}\sum_{p\leq T}\left| \frac{ \mathbf{E}\cos (TU\log p  ) }{\sqrt{p}}\right|\nonumber\\
&:=& I_{1,T}+ I_{2,T}+ I_{3,T}\label{18d-4}
\end{eqnarray}
where we used (\ref{15d-8}). We estimate the three summands. The bound for  $I_{1,T}$ has been obtained in Theorem \ref{t1}.   From this results, we have
\begin{equation*}
I_{1,T}\leq  C \frac{1}{ \log \log T }.
\end{equation*} 
The summand $I_{2,T}$ can be majorized by using Lemma \ref{l2} with $k=1$ and $j=0$.  It holds that
\begin{equation*}
I_{2,T}\leq C \frac{1}{ \sqrt{\log \log T}}.
\end{equation*}
From (\ref{15d-9}), for $T$ large enough, we clearly have 
\begin{equation*}
I_{3,T}\leq C \frac{1}{T\sqrt{\log \log T}} \sum_{p\leq T} \frac{1}{\sqrt{p} \log p}\leq C\frac{1}{\sqrt{T \log  T \log \log T}}
\end{equation*}
due to (\ref{pn1}). These estimates (and similar estimates for the imaginary part) leads to the following result.

\begin{theorem}\label{t3} With $\zeta, U$ as in (\ref{zeta1}), (\ref{u}) respectively, and with $T$ large enough,
\begin{equation}
\label{4i-2} d_{W} \left( \frac{\log \left| \zeta \left( \frac{1}{2}+iTU\right) \right|}{\sqrt{\frac{1}{2} \log \log T}} , N\right)\leq C \frac{1}{ \sqrt{\log \log T}}, d_{W} \left( \frac {\arg\log  \zeta \left( \frac{1}{2}+iTU\right) }{\sqrt{\frac{1}{2} \log \log T}} , N\right)\leq C \frac{1}{ \sqrt{\log \log T}} 
\end{equation}
where $N\sim N(0,1)$.
\end{theorem}

\begin{remark}
Theorem 3 improves the error bound obtained in \cite{Se4} or \cite{Wahl}. In these references, the right-hand bound in (\ref{4i-2}) is $C \frac{\log \log \log T}{\sqrt{\log \log T}}$ under the Kolmogorov metric.

\end{remark}

\section{Multidimensional Selberg theorem and the rate of convergence}

In this paragraph we give a multidimensional extension of the Selberg central limit theorem. Concretely, we consider the $d+1$ dimensional random vector
\begin{equation}
\label{vt}
\mathbf{V}_{T}:= \left( \frac{1}{ \sqrt{\frac{1}{2}\log \log T}} \left( \log \zeta \left( \frac{1}{2} +\mathbf{i} P_{i} \right)\right) \right) _{i=0,.., d}
\end{equation}
and we analyze the asymptotic distribution of its real and imaginary part. As mentionned in the introduction, it has been proved in \cite{HNY} that for $P_{i}= U e ^ {(\log T)^ {\lambda _{i}}}$, $i=0,..,d$ with $\lambda _{0} <\ldots < \lambda _{d}$, the vector $\mathbf{V}_{T}$ (\ref{vt}) converges in law, as $T\to \infty$, to a $d+1$-dimensional standard complex Gaussian vector. When the space between the points $P_{i}$ is small then the limit of $\mathbf{V}_{T}$ is a Gaussian vector with correlated components. The result is due to \cite{B}. 

We analyze the error bound in the multidimensional Selberg theorem in both cases: when the distance between $P_{i}$ and $P_{j}$ is " big" or " small". In fact, we first consider the case $P_{i}= T(U+i)$, $i=0,..,d$. This is a natural multi-dimensional extension of the celebrated Selberg' s result and it seems that it has not yet proved in the literature.  Here $P_{i}- P_{j}= (i-j)T$ is big enough to avoid the correlation between the components of the limit. Then we treat the case considered in \cite{B} when the evaluation points are less distant one from each other. 

The basic idea is the same: one approximates $\mathbf{V}_{T}$ by a random vector whose components are Dirichlet series of the form (\ref{dir}). Using the techniques of the Malliavin calculus, we obtain the rate of convergence of this approximation to the Gaussian limit  under the Wasserstein metric. Then we deduce a bound for the Wasserstein distance between $\mathbf{V}_{T}$ and the Gaussian limit.

\subsection{Big shifts: convergence to a standard Gaussian vector}
Let us  first treat  the case of  "big shifts", i.e. the distance between the evaluation points is big enough and the limit distribution in the multidimensional Selbergh theorem is a standard Gaussian vector.

\subsubsection{Error bound for the Dirichlet series}

Consider the $d+1$ dimensional random vector
\begin{equation*}
\mathbf{X}_{T}= ( X_{T} ^{(0)},..., X _{T} ^{(d)})
\end{equation*}
where for every $i=0,.., d$, 
\begin{equation}
\label{xti}
X _{T} ^{(i)} = \sum _{p\leq T} \frac{1}{\sqrt{p} }\left( \cos (T\log p t_{i} ) - \mathbf{E}\cos (T\log p t_{i} )\right)
\end{equation}
where $t_{i} \sim \mathcal{U}  [Ti, T(i+1)]$. 
We analyze the asymptotic limit of $X_{T}$ as $T \to \infty$.  Since $\mathbf{X}_{T}$ is close to the real part of the random vector $\mathbf{V}_{T}$ (\ref{vt}), we will then deduce the asymptotic distribution and the error bound for $\mathbf{V}_{T}$.   In order to use the techniques of the Malliavin calculus, we will assume that
$$t_{i}= T(U+i), \mbox{ for every } i=0,...,d$$
where $U$ is given by (\ref{u}).  Clearly $t_{i} \sim \mathcal{U}  [Ti, T(i+1)]$.

The exists a multidimensional version of the Stein-Malliavin inequality presented in Theorem \ref{t}, see \cite{NPbook}, \cite{NPR}. This bound is given in terms of the Wasserstein distance. Namely, if $F=(F_{0}, ..., F_{d}) $ is a random vector with components in $\mathbb{D} ^{1,4}$ and $ N(0,\Lambda)$ denotes the $d+1$ dimensional Gaussian distribution with covariance matrix $\Lambda =(c_{i,j})_{i,j=0,..,d}$, then 
\begin{equation}
\label{sm}
d_{W} \left(F,N(0,\Lambda)\right) \leq C \sum_{i=0} ^{d} \mathbf{E} \left|  2c_{i,j}- \langle DF_{i}, D(-L) ^{-1} F_{j}\rangle -\langle DF_{j}, D(-L) ^{-1} F_{i}\rangle  \right|.
\end{equation}
Actually,  the bound presented in \cite{NPbook} or \cite{NPR} is slightly different (and uses the $L^ {2}$-norm on the right-hand side of the inequality (\ref{sm}) )but it is easy to obtain (\ref{sm}) by similar arguments. Recall that the Wasserstein distance between the laws of two $\mathbb{R} ^{d}$ - valued random variables $F,G$ is defined by
\begin{equation}
\label{w}
d_{W} (F,G)= \sup _{h\in \mathcal{A}} \left| \mathbf{E} h(F) - \mathbf{E} h(G) \right|
\end{equation}
where we denote by $\mathcal{A}$ the class of all functions $h: \mathbb{R} ^{d}\to \mathbb{R}$ such that $\Vert h\Vert _{Lip}\leq 1,$ where $\Vert h\Vert _{Lip}= \sup_{x,y\in \mathbb{R} ^{d}, x\not=y} \frac{ \vert h(x)-h(y)\vert}{\Vert x-y\Vert },$
with the  Euclidean norm $\Vert \cdot\Vert $ in $\mathbb{R} ^{d}$. 

In order to apply (\ref{sm}),  we need to calculate 
 $$\langle DX_{T} ^{(i)}, D(-L) ^{-1} X_{T} ^{(j)}\rangle+ \langle DX_{T} ^{(j)}, D(-L) ^{-1} X_{T} ^{(i)}\rangle$$
for every $i,j=0,..,d$. Using the series expansion of the cosinus function, we have
\begin{equation*}
X_{T} ^{(i)} =\sum_{p\leq T} \frac{1}{\sqrt{p}} \sum _{k\geq 0} \frac{ (-1) ^{k} }{(2k)! } (T\log p) ^{2k} G_{2k} ^{(i) }
\end{equation*}
with the notation, for $i=0,..,d$ and for $k>0$
\begin{equation}
\label{gki}
 G_{2k} ^{(i)} = (U+i) ^{2k}- \mathbf{E} (U+i) ^{2k}.
\end{equation}
Then
\begin{eqnarray}
&&\langle DX_{T} ^{(i)}, D(-L) ^{-1} X_{T} ^{(j)}\rangle+ \langle DX_{T} ^{(j)}, D(-L) ^{-1} X_{T} ^{(i)}\rangle =\sum _{p_{1}, p_{2}\leq T} \frac{1}{\sqrt{p_{1}p_{2}}} \sum_{k,l\geq 0} \frac{(-1) ^{k+l} }{ (2k)! (2l)! }\nonumber \\
&&\times  (T\log p_{1}) ^{2k} (T\log p_{2}) ^{2l} \left[ \langle  DG_{2k} ^{(i) } , D(-L) ^{-1} G_{2l} ^{(j)} \rangle + \langle DG_{2l} ^{(j) } , D(-L) ^{-1} G_{2k} ^{(i)} \rangle \right].\label{18d-1}
\end{eqnarray}
The next step is to calculate $\langle  DG_{2k} ^{(i) } , D(-L) ^{-1} G_{2l} ^{(j)} \rangle +  DG_{2l} ^{(j) } , D(-L) ^{-1} G_{2k} ^{(i)} \rangle $. This will be done in the following lemma, based on Lemma \ref{l1}.

\begin{lemma}\label{l3}
For every $k,l>0$ and $i,j=0,..,d$ we have

\begin{eqnarray}
&&\langle DG_{2k} ^{(i) } , D(-L) ^{-1} G_{2l} ^{(j)} \rangle +\langle   DG_{2l} ^{(j) } , D(-L) ^{-1} G_{2k} ^{(i)} \rangle\nonumber\\
&=& \frac{2k}{2l+1}(U+i) ^ {2k-1} U \left[ (j+1) ^{2l+1} -j^{2l+1} -\frac{(U+j) ^{2l+1} }{U}+ \frac{j^{2l+1}}{U} \right] \nonumber\\
&&+  \frac{2l}{2k+1}(U+j) ^ {2l-1} U   \left[ (i+1) ^{2k+1} -i^{2k+1} -\frac{(U+i) ^{2k+1} }{U}+ \frac{i^{2k+1}}{U} \right] \label{24d-5}
\end{eqnarray}
where $G_{2k}^ {(i)}$ is given by (\ref{gki}).
\end{lemma} 
{\bf Proof: } Using the Newton formula $G_{2k} ^{(i)}= \sum _{s=0} ^{2k} C_{2k}^{s}G_{s}  i ^{2k-s} $
and 
\begin{eqnarray*}
&&\langle DG_{2k} ^{(i) } , D(-L) ^{-1} G_{2l} ^{(j)} \rangle +  DG_{2l} ^{(j) } , D(-L) ^{-1} G_{2k} ^{(i)} \rangle\\
&=& \sum _{s=0} ^{2k} \sum _{t=0} ^{2l}C_{2k}^{s}C_{2l} ^{t} i ^{2k-s} j^{2l-t} \left[ \langle DG_{s}, D(-L) ^{-1} G_{t} \rangle + \langle  D G_{t} , D(-L) ^{-1} G_{s} \rangle \right]
\end{eqnarray*}
with $G_{s}$ from (\ref{gk}). By Lemma \ref{l1},
$$\langle DG_{s}, D(-L) ^{-1} G_{t} \rangle + \langle  D G_{t} , D(-L) ^{-1} G_{s} \rangle  =\frac{s}{t+1}U^{s} (1-U ^{t}) + \frac{t}{s+1} U ^{t} (1-U ^{s}). $$
Therefore,
\begin{eqnarray}
&&\langle DG_{2k} ^{(i) } , D(-L) ^{-1} G_{2l} ^{(j)} \rangle +  DG_{2l} ^{(j) } , D(-L) ^{-1} G_{2k} ^{(i)} \rangle\nonumber\\
&=& \sum _{s=0} ^{2k} \sum _{t=0} ^{2l}C_{2k}^{s}C_{2l} ^{t} i ^{2k-s} j^{2l-t}\frac{s}{t+1}U^{s} (1-U ^{t}) + \frac{t}{s+1} U ^{t} (1-U ^{s})\nonumber\\
&=& \sum _{s=0} ^{2k} C_{2k}^{s}i ^{2k-s} sU^{s} \sum _{t=0} ^{2l} C_{2l} ^{t} j^{2l-t} \frac{1-U^{t}}{t+1}+ \sum _{s=0} ^{2k} C_{2k}^{s}i ^{2k-s} \frac{1-U^{s}}{s+1} \sum _{t=0} ^{2l} C_{2l} ^{t} j^{2l-t} tU^{t}.\label{24d-6}
\end{eqnarray}
Now we calculate the sums after $s$ and $t$.  Notice that

$$\sum _{s=0} ^{2k} C_{2k}^{s}i ^{2k-s} sU^{s}= 2k (U+i) ^{2k-1} U$$
and
$$\sum _{t=0} ^{2l} C_{2l} ^{t} j^{2l-t} \frac{1-U^{t}}{t+1} = \frac{1}{2l+1} \left[ (j+1) ^{2l+1} -j^{2l+1} -\frac{(U+j) ^{2l+1} }{U}+ \frac{j^{2l+1}}{U} \right]. $$

From the above two identities and (\ref{24d-6}), we deduce the conclusion (\ref{24d-5}).\qed

\begin{remark}
For $i=j=0$, we retrieve the result in Lemma \ref{l1}. 
\end{remark}

We are now in position to compute the terms involving Malliavin operators that appear in the right-hand side of (\ref{sm}).

\begin{lemma}\label{l5} For $i,j=0,..,d$, let $X_{T} ^ {(i)}$ be defined by (\ref{xti}). Then
\begin{eqnarray*}
&&\langle DX_{T} ^{(i)}, D(-L) ^{-1} X_{T} ^{(j)}\rangle+ \langle DX_{T} ^{(j)}, D(-L) ^{-1} X_{T} ^{(i)}\rangle  = \sum _{p_{1}, p_{2}\leq T} \frac{1}{\sqrt{p_{1} p_{2}} }\frac{\log p_{2}} {\log p_1} \\
&&\left[ \sin ((U+i) T \log p_{1})\sin ((U+j)T \log p_{2})+\sin ((U+j) T \log p_{1})\sin ((U+i)T \log p_{2})\right]  +R_{T}^ {(i,j)}
\end{eqnarray*}
with
\begin{eqnarray*}
R_{T}^ {(i,j)}&=& \sum _{p_{1}, p_{2}\leq T} \frac{1}{\sqrt{p_{1} p_{2}} }\frac{\log p_{2}} {\log p_1} \\
&&\left[ -\frac{U}{i+1} \sin ((U+j) T \log  p_{2})\sin( (i+1) T\log p_{1}) - \frac{U}{j+1} \sin ((U+i) T \log  p_{2})\sin( (j+1) T\log p_{1}) \right.\\
&&\left. 
+\frac{U}{i} \sin ((U+j) T \log p_{2}) \sin (iT \log p_{1}) 1_{i\not =0} +\frac{U}{j} \sin ((U+i) T \log p_{2}) \sin (jT \log p_{1}) 1_{j\not =0}\right.\\
&&\left.   -\sin ((U+j)T \log p_{2}) \sin (iT \log p_{1}) -\sin ((U+i)T \log p_{2}) \sin (jT \log p_{1}) \right].
\end{eqnarray*}

\end{lemma}
{\bf Proof: } By Lemma \ref{l2} and (\ref{18d-1}),
\begin{eqnarray}
&&\langle DX_{T} ^{(i)}, D(-L) ^{-1} X_{T} ^{(j)}\rangle+ \langle DX_{T} ^{(j)}, D(-L) ^{-1} X_{T} ^{(i)}\rangle \nonumber \\
&=& \sum _{p_{1}, p_{2}\leq T} \frac{1}{\sqrt{p_{1}p_{2}}}\sum_{k,l\geq 0} \frac{(-1) ^{k+l} }{ (2k)! (2l)! } (T\log p_{1}) ^{2k} (T\log p_{2}) ^{2l} \nonumber\\
&&\left( \frac{2k}{2l+1}(U+i) ^ {2k-1} U \left[ (j+1) ^{2l+1} -j^{2l+1} -\frac{(U+j) ^{2l+1} }{U}+ \frac{j^{2l+1}}{U} \right]\right. \nonumber\\
&&\left. +  \frac{2l}{2k+1}(U+j) ^ {2l-1} U   \left[ (i+1) ^{2k+1} -i^{2k+1} -\frac{(U+i) ^{2k+1} }{U}+ \frac{i^{2k+1}}{U} \right]\right) \nonumber\\
&=&\sum _{p_{1}, p_{2}\leq T} \frac{1}{\sqrt{p_{1}p_{2}}} \sum _{k\geq 0} \frac{(-1)^ {k}}{(2k+1) !} (T\log p_{1}) ^ {2k}  \left[ (i+1) ^{2k+1} -i^{2k+1} -\frac{(U+i) ^{2k+1} }{U}+ \frac{i^{2k+1}}{U} \right]\nonumber\\
&&\times U\sum _{l\geq 1} \frac{(-1) ^ {l} }{(2l-1) ! } (U+j) ^ {2l-1} (T\log p_{2} ) ^ {2l} \nonumber\\
&+&\sum _{p_{1}, p_{2}\leq T} \frac{1}{\sqrt{p_{1}p_{2}}}  \sum _{l\geq 0} \frac{(-1)^ {l}}{(2l+1) !} (T\log p_{1}) ^ {2l}  \left[ (j+1) ^{2l+1} -i^{2l+1} -\frac{(U+j) ^{2l+1} }{U}+ \frac{j^{2l+1}}{U} \right]\nonumber\\
&&\times U\sum _{k\geq 1} \frac{(-1) ^ {k} }{(2k-1) ! } (U+i) ^ {2k-1} (T\log p_{2} ) ^ {2k}. \label{24d-8}
\end{eqnarray}
Next, we calculate the above  sums over $k$ and $l$. We have
\begin{equation}
\label{24d-9}
\sum _{l\geq 1} \frac{(-1) ^ {l} }{(2l-1) ! } (U+j) ^ {2l-1} (T\log p_{2} ) ^ {2l} =-U T\log p_{2} \sin ( (U+j) T \log p_{2})
\end{equation}
and

\begin{eqnarray}
&&\sum_{k\geq 0}\frac{(-1)^ {k}}{(2k+1) !} (T\log p_{1}) ^ {2k}  \left[ (i+1) ^{2k+1} -i^{2k+1} -\frac{(U+i) ^{2k+1} }{U}+ \frac{i^{2k+1}}{U} \right] \nonumber\\
&=&
\frac{1}{T\log p_{1}} \left[ \frac{1}{i+1} \sin( T(i+1) \log p_{1} ) )\right. \nonumber\\
&&\left.  - \frac{1}{i} \sin (Ti\log p_{1} ) -\frac{1}{U} \sin(T(U+i)\log p_{1} ) + \frac{1}{U} \sin (Ti \log p_{1}  )\right]. \label{24d-10}
\end{eqnarray}
By plugging relations (\ref{24d-9}) and (\ref{24d-10}) into (\ref{24d-8}), we get the conclusion.\qed

\vskip0.2cm

We measure now the Wasserstein distance between the (renormalized) sequence $\mathbf{X}_{T}$ and the standard $d+1$ dimensional Gaussian distribution (denoted $N(0, I_{d+1})$ in the sequel).

\begin{prop} \label{p3}Let $ X_{T} ^ {(i)}$ be given by (\ref{xti}) for $i=0,..,d$ and let
Let $$ \mathbf{F}_{T}=  \sqrt{\frac{1}{\frac{1}{2}\log \log T}}\mathbf{ X}_{T}=   \sqrt{\frac{1}{\frac{1}{2}\log \log T}} (X _{T} ^ {(0)}, X_{T} ^ {(1)},..., X _{T} ^ {(d)}). $$
Then, for large $T$
\begin{equation}\label{18d-3}
d _{W} (\mathbf{F}_{T}, N(0, I_{d+1})) \leq C \frac{1}{\log \log T}. 
\end{equation}
\end{prop}
{\bf Proof: } Let $F ^{(j)}_{T}, j=0,..,d$ be the components of the vector $\mathbf{F}_{T}$. Using the Stein-Malliavin bound (\ref{sm}), we have
\begin{eqnarray*}
d_{W} (\mathbf{F}_{T}, N(0, I_{d+1})) &\leq &  C \left[ \sum _{i=0} ^ {d+1} \mathbf{E}\left| 1- \langle DF_{T} ^ {(i) } , D(-L)^ {-1} F _{T} ^ {(i)}\rangle \right| \right. \\
&+&  \left.   \sum _{i=0} ^ {d+1} \mathbf{E} \left| \langle DF_{T} ^ {(i) } , D(-L)^ {-1} F _{T} ^ {(j)}\rangle +\langle DF_{T} ^ {(j) } , D(-L)^ {-1} F _{T} ^ {(i)}\rangle \right|\right].
\end{eqnarray*}
The main contribution will come from the diagonal term. By Lemma \ref{l5} (with $R^ {(i,j})_{T}$ as in the statement of Lemma \ref{l5}) ,
\begin{eqnarray*}
&&\left|  \langle DF_{T} ^ {(i) } , D(-L)^ {-1} F _{T} ^ {(i)}\rangle -1\right|\\
&=& \frac{1}{\frac{1}{2} \log \log T} \sum _{p\leq T} \frac{1}{p} \sin ^ {2}((U+i) T \log p)\\
&&+ \frac{1}{\frac{1}{2} \log \log T} \mathbf{E}\sum_{p_{1}, p_{2}\leq T; p_{1}\not=p_{2}} \frac{1}{\sqrt{p_{1} p_{2}} }\frac{\log p_{2}} {\log p_1} \sin ((U+i) T \log p_{1})\sin ((U+j)T \log p_{2}) \\
&&+ \frac{1}{\frac{1}{2} \log \log T} R_{T} ^ {(i,i)} -1\\
&=& \left[ \frac{1}{\frac{1}{2} \log \log T} \sum_{p } \frac{1}{2p}  -1\right]-\frac{1}{\frac{1}{2} \log \log T} \sum _{p} \frac{1}{2p} \cos  (2(U+i) T \log p)\\
&&+\frac{1}{\frac{1}{2} \log \log T}  \sum_{p_{1}, p_{2}\leq T;p_{1}\not=p_{2}} \frac{1}{\sqrt{p_{1} p_{2}} }\frac{\log p_{2}} {\log p_1} \sin ((U+i) T \log p_{1})\sin ((U+j)T \log p_{2}) \\
&&+ \frac{1}{\frac{1}{2} \log \log T} R_{T} ^ {(i,i)} \\
&:=&  \left[ \frac{1}{\frac{1}{2} \log \log T} \sum_{p } \frac{1}{2p}  -1\right] + R_{T}
\end{eqnarray*}
where we included in the rest term $R_{T}$ the summands $-\frac{1}{\frac{1}{2} \log \log T} \sum _{p} \frac{1}{2p} \cos  (2(U+i) T \log p)$, $\frac{1}{\frac{1}{2} \log \log T}  \sum_{p_{1}\not=p_{2}} \frac{1}{\sqrt{p_{1} p_{2}} }\frac{\log p_{2}} {\log p_1} \sin ((U+i) T \log p_{1})\sin ((U+j)T \log p_{2}) $ and $\frac{1}{\frac{1}{2} \log \log T} R_{T} ^ {(i,i)} $.

\b The first term above is the one which gives the bound (\ref{18d-3}). Indeed, by (\ref{pn2}), 
$$\left|  \frac{1}{\frac{1}{2} \log \log T} \sum_{p\leq T } \frac{1}{2p}  -1\right| \leq C\frac{1}{\log \log T}.$$
On the other hand, using only the bounds (\ref{15d-9}), (\ref{15d-6}) and (\ref{18d-2}), we can easily proof that
$$ \mathbf{E}\vert R_{T}\vert  \leq C\frac{1}{\log \log T}.$$
(actually, a better estimated is possible but not necessarily for our purpose). Again by (\ref{15d-9}), (\ref{15d-6}) and (\ref{18d-2}), 
$$\mathbf{E} \left| \langle DF_{T} ^ {(i) } , D(-L)^ {-1} F _{T} ^ {(j)}\rangle +\langle DF_{T} ^ {(j) } , D(-L)^ {-1} F _{T} ^ {(i)}\rangle \right|\leq  C\frac{1}{\log \log T}$$
for every $i\not=j$. The two above estimates lead to (\ref{18d-3}). \qed 

\subsubsection{Error bound for the multidimensional Selberg theorem}

We regard now the asymptotic behavior of the vector 
$$\left(  \log \left| \zeta (\frac{1}{2}+ \mathbf{i} P_{i} )\right| \right)_{i=0,..,d}$$
with the evaluation points $P_{i}= T(U+i), i=0,..,d$. We show that, after normalization, it also converges to a standard Gaussian vector. Recall that $U$ denote a standard uniform random variable defined by (\ref{u}).

We have 
\begin{theorem}\label{t4}
Let  
\begin{equation}
\label{xxt}
{\bf  \mathcal{X}} _{T}= \frac{1}{\sqrt{\frac{1}{2}\log \log T}}\left( \log \vert \zeta (\frac{1}{2}+\mathbf{i}UT)\vert,  \log \vert \zeta (\frac{1}{2}+\mathbf{i}(U+1)T)\vert, \ldots, \log \vert \zeta (\frac{1}{2}+\mathbf{i}(U+d)T)\vert \right)
\end{equation}
with $U$ from (\ref{u}). Then for large $T$,
$$d_{W}({\bf \mathcal{ X}} _{T}, N(0, I_{d+1})\leq C \frac{1}{\sqrt{\log \log T}}.$$
\end{theorem}
{\bf Proof: } As in (\ref{18d-4}), using the triangle inequality for the Wasserstein distance and the fact that for any $\mathbb{R}^ {d}$- valued random variables $F,G$ we have 
$$d_{W}(F,G)\leq \mathbf{E} \vert F-G\vert _{1} \leq\left( \mathbf{E} \vert F-G\vert _{2} \right) ^ {\frac{1}{2}}$$
where $\vert x\vert _{1} = \vert x_{1}\vert+....\vert x_{d}\vert$ and $\vert x\vert _{2}^ {2} =\vert x_{1}\vert ^ {2} +...+ \vert x_{d}\vert  ^ {2}$ if $x=(x_{1},.., x_{d})$, we can write
\begin{eqnarray*}
d_{W}({\bf \mathcal{X}} _{T}, N(0, I_{d+1}) &\leq & d_{W} (\mathbf{F}_{T}, N(0, I_{d+1})) \\
&+& \mathbf{E} \vert {\bf \mathcal {X}} _{T}- \mathbf{F}_{T} \vert _{1} + \frac{1}{\sqrt{\frac{1}{2}\log \log T}}\sum_{p\leq T} \frac{1}{\sqrt{p}}\left| \left( \mathbf{E} (\cos (T\log p (U+i) ))\right) _{i=0,..,d} \right| _{1}
\end{eqnarray*}
and by Proposition \ref{p3}
$$ d_{W} (\mathbf{F}_{T}, N(0, I_{d+1})) \leq C\frac{1}{\log \log T} $$
while the inequality (\ref{15d-7}) in Lemma \ref{l2} implies 
$$ \mathbf{E} \vert \mathcal {X} _{T}- \mathbf{F}_{T} \vert _{1}  \leq C \frac{1}{\sqrt{\log \log T}}.$$
Finally, as in the proof of Theorem \ref{t1} (using (\ref{15d-9}) and (\ref{pn2}))
$$\frac{1}{\sqrt{\frac{1}{2}\log \log T}}\sum_{p\leq T} \frac{1}{\sqrt{p}}\left| \left( \mathbf{E} (\cos (T\log p (U+i) ))\right) _{i=0,..,d} \right| _{1} \leq C\frac{1}{\sqrt{\log \log T}} .$$

\qed

Exactly the same analysis can be done for the imaginary part of $\log \zeta $ on the critical line.  For $T>0$, let  us define 
$$\mathbf{Y}_{T}= ( Y_{T} ^ {(0)},..., Y_{T}^ {(d)})$$ 
where, for every $i=0,..,d$, 
$$Y_{T} ^ {(i)}= \sum _{p\leq T} \frac{1}{\sqrt{p}}\left[ \sin (T\log p (U+i)) -\mathbf{E} \sin (T\log p (U+i)) \right]. $$
As in the proof of Lemma \ref{l5}, we can see that for every $i,j=0,..,d$,
\begin{eqnarray*}
&&\langle DY_{T} ^{(i)}, D(-L) ^{-1} Y_{T} ^{(j)}\rangle+ \langle DY_{T} ^{(j)}, D(-L) ^{-1} Y_{T} ^{(i)}\rangle  = \sum _{p_{1}, p_{2}\leq T} \frac{1}{\sqrt{p_{1} p_{2}} }\frac{\log p_{2}} {\log p_1}\\
&&\left[ \cos ((U+i) T \log p_{1})\cos ((U+j)T \log p_{2})  +\cos ((U+i) T \log p_{2})\cos ((U+j)T \log p_{1})   \right] \\
&&+R_{T}^ {(i,j)}
\end{eqnarray*}
with
\begin{eqnarray*}
R_{T}^ {(i,j)}&=& \sum _{p_{1}, p_{2}\leq T} \frac{1}{\sqrt{p_{1} p_{2}} }\frac{\log p_{2}} {\log p_1}\\
&&\left[ -\cos(T(U+j)\log p_{2}) \cos (Ti\log p_{1} )  -\cos(T(U+i)\log p_{1}) \cos (Tj\log p_{1} )\right. \\
&& \left. -U  \cos(T(U+j)\log p_{2}) \cos (T(i+1) \log p_{1} )-U \cos(T(U+i)\log p_{1}) \cos (T (j+1)\log p_{1})\right.\\
&&\left.+ U \cos(T(U+j)\log p_{2}) \cos (Ti\log p_{1} )+U \cos(T(U+i)\log p_{2}) \cos (Tj\log p_{1} )\right].
\end{eqnarray*}
Using the proof of Proposition \ref{p3}, with 
$$\mathbf{G}_{T}:= \frac{1}{\sqrt{\frac{1}{2}\log \log T}} \mathbf{Y}_{T},$$
 we can show that 
$$d_{W} (\mathbf{G}_{T}, N(0, I_{d+1}))\leq C\frac{1}{\log \log T. }$$ We will then obtain a similar result to Theorem \ref{t4}.

\begin{theorem}
With $U$ from (\ref{u}) and $\zeta$ from (\ref{zeta1}), let  
\begin{equation}
\label{yyt}
 \mathcal{Y} _{T}= \frac{1}{\sqrt{\frac{1}{2}\log \log T}}\left(\arg  \log  \zeta (\frac{1}{2}+\mathbf{i}UT), \arg \log \zeta (\frac{1}{2}+\mathbf{i}(U+1)T), \ldots, \arg \log  \zeta (\frac{1}{2}+\mathbf{i}(U+d)T) \right).
\end{equation}
Then for large $T$,
$$d_{W}(\mathcal{Y} _{T}, N(0, I_{d+1})\leq C \frac{1}{\sqrt{\log \log T}}.$$
\end{theorem}

\subsection{Small shifts: Convergence to a Gaussian vector with correlated components}

The next step is to analyze the asymptotic behavior of the vector $\mathbf{V}_{T}$ (\ref{t}) when the evaluation points $P_{i}$ are close one to each other. That is, we choose $P_{i} =TU+ f_{T}^ {(i)}$, $i=0,.., d$, with $U$ a standard uniform random variable   defined by (\ref{u}) and $f_{T}^ {(i)}$ are small shifts, meaning that  the difference $f_{T} ^ {(i)}- f_{T} ^ {(j)}$ is small enough for $i\not=j$. This will lead to the appearance of non trivial correlations between the components of the limit of (\ref{vt}) as $T\to \infty$.

We first look to the asymptotic behavior of the $d+1$-dimensional  Dirichlet series that approximates $\mathbf{V}_{T}$.  That is, we introduce 
\begin{equation*}
\mathbf{Z}_{T}= \left( Z_{T} ^ {(0)}, Z _{T} ^ {(1)}, \ldots , Z_{T} ^ {(d) }\right)
\end{equation*}
with

\begin{equation}\label{zti}
 Z_{T}^ {(i)}=\sum_{p\leq T} \cos\left(   (TU + f_{T} ^ {(i)})\log p\right), \hskip0.3cm i=0,..,d
\end{equation}
We will give the rate of convergence of $\mathbf{Z}_{T}$ to the $d+1$-dimensional Gaussian law with suitable covariance matrix. As before, the proof will be based on the inequality (\ref{sm}) and it  means that we need to compute $\langle DZ _{T} ^ {(i)} , D(-L) ^ {-1}Z _{T} ^ {(j) }\rangle $ for every $i,j=0,..,d$. We have

\begin{lemma}
For $i=1,,,,d$, let   $f_{T} ^ {(i)} $ be a deterministic function and let $Z_{T} ^ {(i)} $ be given by (\ref{zti}). Then
\begin{eqnarray*}
&&\langle DZ _{T} ^ {(i)} , D(-L) ^ {-1}Z _{T} ^ {(j) }\rangle +\langle DZ _{T} ^ {(j)} , D(-L) ^ {-1}Z _{T} ^ {(i) }\rangle \\
&=&\sum _{p_{1}, p_{2}\leq T} \frac{1}{\sqrt{p_{1}p_{2}}} \frac{\log p_{2}}{\log p_{1}}\\
&&\left[   \sin\left( \log p_{1} (TU+ f_{T} ^ {(i) }) \right) \sin \left( \log p_{2}(TU+ f_{T} ^ {(j) }) \right) + \sin\left( \log p_{1} (TU+ f_{T} ^ {(j) }) \right) \sin \left( \log p_{2}(TU+ f_{T} ^ {(i) }) \right) \right]\\
&&+ R^ {(i,j)}_{T}
\end{eqnarray*}
with
\begin{eqnarray*}
&&R^ {(i,j)}_{T}=\sum _{p_{1}, p_{2}\leq T} \frac{1}{\sqrt{p_{1}p_{2}}} \frac{\log p_{2}}{\log p_{1}}\\
&&\left[  (U-1)\sin\left( (TU+ f_{T} ^ {(j) }) \log p_{1}\right) \sin \left( f_{T} ^ {(i)}\log p_{2}\right)+ (U-1)\sin\left( (TU+ f_{T} ^ {(i) })\log p_{1} \right) \sin\left( f_{T} ^ {(j)}\log p_{2}\right)\right. \\
&&\left. -U\sin\left( (TU+ f_{T} ^ {(j) }) \log p_{1}\right) \sin \left((T+f_{T} ^ {(i)}) \log p_{2}\right) -U\sin\left(  (TU+ f_{T} ^ {(i) }) \log p_{1} \right)\sin \left((T+f_{T} ^ {(j)})\log p_{2} \right)\right].
\end{eqnarray*}

\end{lemma}
{\bf Proof: }With arguments previously used, we can write

\begin{eqnarray*}
&&\langle DZ _{T} ^ {(i)} , D(-L) ^ {-1}Z _{T} ^ {(j) }\rangle +\langle DZ _{T} ^ {(j)} , D(-L) ^ {-1}Z _{T} ^ {(i) }\rangle 
=\sum _{p_{1}, p_{2}} \frac{1}{\sqrt{p_{1}p_{2}}}\sum_{k,l\geq 0} \frac{(-1) ^{k+l} }{ (2k)! (2l)! } (\log p_{1}) ^{2k} (\log p_{2}) ^{2l} \nonumber\\
&&\times \left[  DH_{2k} ^{(i) } , D(-L) ^{-1} H_{2l} ^{(j)} \rangle +  \langle DH_{2l} ^{(j) } , D(-L) ^{-1} H_{2k} ^{(i)} \rangle \right]
\end{eqnarray*}
where we used the notation
\begin{equation}
\label{hki}
H_{2k} ^ {(i)}= (TU+ f_{T} ^ {(i)} ) ^ {2k}-\mathbf{E}(TU+ f_{T} ^ {(i)} ) ^ {2k}.
\end{equation}
The scalar product above will be computed as in the proof of Lemma \ref{l3}. We get
\begin{eqnarray*}
&&\langle DH_{2k} ^{(i) } , D(-L) ^{-1} H_{2l} ^{(j)} \rangle +  \langle DH_{2l} ^{(j) } , D(-L) ^{-1} H_{2k} ^{(i)} \rangle\\
&=&\frac{2k}{2l+1} U (TU+ f_{T } ^ {(i)} ) ^ {2k-1} \left[ (T+f_{T}^ {(j)}) ^ {2l+1} -( f_{T} ^ {(j)}) ^ {2l+1}-\frac{1}{U}  (T+f_{T}^ {(j)}) ^ {2l+1} + \frac{1}{U} ( f_{T} ^ {(j)}) ^ {2l+1}\right]\\
&+& \frac{2l}{2k+1} U (TU+ f_{T } ^ {(j)} ) ^ {2l-1} \left[ (T+f_{T}^ {(i)}) ^ {2k+1} -( f_{T} ^ {(i)}) ^ {2k+1}-\frac{1}{U}  (TU+f_{T}^ {(i)}) ^ {2k+1} + \frac{1}{U} ( f_{T} ^ {(i)}) ^ {2k+1}\right].
\end{eqnarray*}
When we compute the above sums after $k$ and $l$, we obtain 
$$\sum_{l\geq 1} U \frac{(-1) ^ {l} }{(2l-1) ! }(\log p_{2}) ^ {2l} (TU+ f_{T } ^ {(j)} ) ^ {2l-1} =-U\log p_{2} \sin \left( \log p_{2}(TU+ f_{T} ^{(j)})\right)$$
\begin{eqnarray*}
&&\sum _{k\geq 0}  \frac{(-1) ^ {k}}{(2k+1) !} (\log p_{1}) ^ {2k} \left[ (T+f_{T}^ {(i)}) ^ {2k+1} -( f_{T} ^ {(i)}) ^ {2k+1}-\frac{1}{U}  (TU+f_{T}^ {(i)}) ^ {2k+1} + \frac{1}{U} ( f_{T} ^ {(i)}) ^ {2k+1}\right]\\
&=& \frac{1}{\log p_{1}} \left[ \sin \left( \log p_{1}(T + f_{T} ^ {(i)}) \right)- \sin\left( \log p_{1} (f_{T} ^ {(i)})\right) - \frac{1}{U} \sin\left( \log p_{1} (TU+ f_{T} ^ {(i) }) \right)+ \frac{1}{U} \sin \left(\log p_{1}(f_{T} ^ {(i)})\right)\right].
\end{eqnarray*}
The conclusion follows easily. \qed

Before stating the main results of this section, let us recall the following technical result due to \cite{B}, which plays a key role.

\begin{lemma}\label{l6}
Let $ (\Delta _{T})_{T} $ be bounded and positive such that $\frac{\log \Delta _{T}}{\log \log T} \to _{T}c\in [0, \infty].$ Then $\frac{1}{\log \log T} \sum _{p\leq T} \frac{\cos (\log p \log \Delta  _{T})}{p} \to _{T} c\wedge 1 $
and
\begin{equation}
\label{4i-3}
\left| \frac{1}{\log \log T} \sum _{p\leq T} \frac{\cos (\log p \log \Delta  _{T})}{p} - c\wedge 1\right| \leq C \frac{1}{\log \log T}.
\end{equation}
\end{lemma}
The bound (\ref{4i-3}) is not explicitly stated in \cite{B}, but its proof is an easy consequence of the proof of Lemma 3.4 in \cite{B}.

\vskip0.3cm 
The main result of this paragraph states as follows.

\begin{prop}\label{p4}
Assume $0\leq f_{T} ^ {(0)}<f_{T}^ {(1)} <\ldots f_{T} ^  {(d)} <C <\infty$.  For every $i,j=0,..,d$ with $i\not=j$ suppose that
\begin{equation}\label{28d-1}
\frac{\log \vert f_{T} ^ {(i)} -f_{T} ^ {(j)}\vert }{\log \log T } \to a_{i,j}\in [0, \infty].
\end{equation}
Define
$$\mathbf{A}_{T}:= \frac{1}{\sqrt{\frac{1}{2}\log \log T}} \mathbf{Z}_{T}$$
where $\mathbf{Z}_{T}$ is the vector with components (\ref{zti}). Then $\mathbf{A} _{T}$ converges in distribution, as $T \to \infty$, to a centered Gaussian vector with covariance matrix $\Lambda= (c_{i,j})_{i,j=0,..,d}$ with $c_{i,j}=a_{i,j}\wedge 1$. Moreover, for $T$ large
$$d_{W} ( \mathbf{A}_{T}, N(0,\Lambda))\leq C \frac{1}{\log \log T}.$$
\end{prop}
{\bf Proof: } From (\ref{sm})
\begin{eqnarray*}
&&d_{W} ( A_{T}, N(0,\Lambda ))\\
&\leq & C \sum_{i=0} ^ {d} \mathbf{E} \left|1 -\frac{1}{\frac{1}{2} \log \log T}\sum _{p } \frac{1}{p}  \sin ^ {2} \left(\log p(TU+ f_{T} ^ {(i)})\right)\right| \\
&&+ C \sum _{i\not=j} \mathbf{E} \left| c_{i,j}-\frac{1}{\frac{1}{2} \log \log T} \sum _{p} \frac{1}{p} \sin\left(\log p (TU+ f_{T} ^ {(i)})\right) \sin \left( \log p (TU+ f_{T} ^ {(j)})\right) \right|\\
&& +r_{T}
\end{eqnarray*}
with
\begin{eqnarray*}
r_{T}&=& \frac{1}{\log \log T} \sum_{i,j}\mathbf{E}\vert  R_{T} ^ {(i,j)}\vert \\
&&+ \frac{1}{\log \log T}\sum_{i=0}^ {d} \sum _{p_{1}\not=p_{2} } \frac{1}{\sqrt{p_{1}p_{2}}}\frac{\log p_{1}}{\log p_{2}}
\left| \mathbf{E} \sin\left( \log p_{1} (TU+ f_{T} ^ {(i) }) \right) \sin \left( \log p_{2}(TU+ f_{T} ^ {(i) }) \right) \right|\\
&&+\frac{1}{\log \log T}\sum_{i,j=0; i\not=j}^ {d} \sum _{p_{1}\not=p_{2} } \frac{1}{\sqrt{p_{1}p_{2}}}\frac{\log p_{1}}{\log p_{2}} \left| \mathbf{E}\sin\left( \log p_{1} (TU+ f_{T} ^ {(i) }) \right) \sin \left( \log p_{2}(TU+ f_{T} ^ {(j) }) \right)\right| .
\end{eqnarray*}
Hence, with the trigonometric indentity $\sin ^ {2} (x)= \frac{1- \cos (2x)}{2}$

\begin{eqnarray}
d_{W} (\mathbf{ A}_{T}, N(0,\Lambda)) &\leq & C  \left|1 -\frac{1}{\log \log T} \sum _{p\leq T } \frac{1}{p} \right| \nonumber\\
&&+  C \sum_{i,j=0; i\not=j} ^ {d}\left| c_{i,j}-\frac{1}{\log \log T}  \sum _{p\leq T} \frac{1}{p} \cos \left( \log p (f_{T} ^ {(i)} - f_{T} ^ {(j)} )\right)\right|\nonumber \\
&&+ r_{T,2}+ r_{T}\label{28d-2}
\end{eqnarray}
with
\begin{eqnarray*}
r_{T,2}&=& C \frac{1}{\log \log T} \sum_{i=0} ^ {d} \sum_{p\leq T} \frac{1}{p}\left| \mathbf{E} \cos \left( 2\log p (TU+ f_{T} ^ {(i)} )\right) \right| \\
&&+   C \frac{1}{\log \log T} \sum_{i,j=0; i\not=j} ^ {d}\sum_{p\leq T} \frac{1}{p}\left| \mathbf{E} \cos \left( \log p( 2TU + f_{T} ^ {(i)}+ f_{T}^ {(j)} )\right)  \right|.
\end{eqnarray*}
By (\ref{pn2}),
\begin{equation}
\label{28d-3}
    \left|1 -\frac{1}{\log \log T} \sum _{p \leq T} \frac{1}{p} \right| \leq C \frac{1}{\log \log T}.
\end{equation}
Using
\begin{eqnarray*}
&&\left| \mathbf{E}\cos (C\log p (TU+ a_{T})) \right| \leq C \frac{1}{T\log p}, \left| \mathbf{E}\sin (C\log p (TU+ a_{T})) \right| \leq C \frac{1}{T\log p} \nonumber\\
&&\left| \mathbf{E}U\sin (C\log p (TU+ a_{T})) \right| \leq C \frac{1}{T\log p}
\end{eqnarray*}
we immediately get
\begin{equation}\label{28d-4}
\mathbf{E} \vert r_{T}\vert \leq  C \frac{1}{\log \log T} \mbox{ and } \mathbf{E} \vert r_{2,T}\vert \leq  C \frac{1}{\log \log T} .
\end{equation}
By applying Lemma \ref{l6} to $\Delta _{T}= f_{T} ^ {(j)}-f_{T}^ {(i)}$ we get 
\begin{equation}
\label{28d-5}\frac{1}{\log \log T} \sum_{i,j=0; i\not=j} ^ {d}\left| c_{i,j}- \sum _{p\leq T} \frac{1}{p} \cos \left( \log p (f_{T} ^ {(i)} - f_{T} ^ {(j)} )\right)\right|\leq  C \frac{1}{\log \log T}.
\end{equation}
By inserting (\ref{28d-3}), (\ref{28d-4}) and (\ref{28d-5}) into (\ref{28d-2}), we obtain the desired conclusion. \qed

Concerning the imaginary part of the vector $\mathbf{V}_{T}$ (\ref{vt}), we let
\begin{equation*}
\mathbf{B}_{T}= \frac{1}{\sqrt{\frac{1}{2}\log \log T}}\left( W_{T} ^ {(0)},\ldots , W_{T} ^ {(d)} \right)
\end{equation*}
with
$$ W^ {(i)}_{T}= \sum _{p} \left[ \sin \left( \log p (TU + f_{T} ^ {(i)})\right) - \mathbf{E} \sin \left( \log p (TU + f_{T} ^ {(i)})\right) \right.]$$
Then we will get for every $i,j=0, ..,d$
\begin{eqnarray*}
&&\langle D W _{T} ^ {(i)}, D(-L) ^ {-1} W_{T} ^ {(j)}\rangle \\
&=&\sum _{p_{1}, p_{2}\leq T} \frac{1}{\sqrt{p_{1}p_{2}}} \frac{\log p_{2}}{\log p_{1}}\\
&&\left[  \cos \left( \log p_{2} (TU+ f_{T} ^ {(j)})\right)  \cos \left( \log p_{1} (TU+ f_{T} ^ {(i)})\right) +\cos \left( \log p_{2} (TU+ f_{T} ^ {(i)})\right)  \cos \left( \log p_{1} (TU+ f_{T} ^ {(j)})\right) \right] \\
&&+R ^ {(i,j)}_{T}
\end{eqnarray*}
where
\begin{eqnarray*}
&&R^ {(i,j)}_{T}=\sum _{p_{1}, p_{2}\leq T} \frac{1}{\sqrt{p_{1}p_{2}}} \frac{\log p_{2}}{\log p_{1}}\\
&&\left[ -U\cos \left( \log p_{2} (TU+ f_{T} ^ {(j)})\right)  \cos \left( \log p_{1} (T+ f_{T} ^ {(i)})\right)- U\cos \left( \log p_{2} (TU+ f_{T} ^ {(i)})\right)  \cos \left( \log p_{1} (T+ f_{T} ^ {(j)})\right)\right. \\
&&\left.  +U\cos \left( \log p_{2} (TU+ f_{T} ^ {(j)})\right)  \cos \left( \log p_{1} f_{T} ^ {(i)}\right)+U\cos \left( \log p_{2} (TU+ f_{T} ^ {(i)})\right)  \cos \left( \log p_{1} f_{T} ^ {(j)}\right)\right].
\end{eqnarray*}
This will lead to the following result:

\begin{prop}\label{p5}
Let the assumptions in Proposition \ref{p4} prevail and let $\mathbf{B}_{T}$ be as above. Then
\begin{equation*}
d_{W}\left( \mathbf{B}_{T}, N(0,\Lambda)\right) \leq C \frac{1}{\log \log T}.
\end{equation*}
\end{prop}

\begin{theorem}\label{t6} Let the assumption in Proposition \ref{p4} prevail. Define
$$\mathcal{A}_{T}= \frac{1}{\sqrt{\frac{1}{2}\log \log T}}\left( \log \left| \zeta \left( \frac{1}{2}+ \mathbf{i}(TU+f_{T}^ {(i)}) \right)\right|\right)_{i=0,..,d}$$
and
$$\mathcal{B}_{T}= \frac{1}{\sqrt{\frac{1}{2}\log \log T}}\left( \arg \log\zeta \left( \frac{1}{2}+ \mathbf{i}(TU+f_{T}^ {(i)}) \right)\right)_{i=0,..,d}.$$
Then
\begin{equation*}
d_{W} (\mathcal{A}_{T}; N(0,\Lambda)) \leq C \frac{1}{\sqrt{\log \log T}} \mbox{ and } d_{W} (\mathcal{B}_{T}; N(0,\Lambda)) \leq C \frac{1}{\sqrt{\log \log T}} .
\end{equation*}
\end{theorem}
{\bf Proof: } The conclusion follows from Propositions \ref{p4} and \ref{p5} and from the fact that the conclusion of Lemma \ref{l2} is true if we replace $T(U+i)$ by $TU+ f_{T} ^ {(i)}$ (see Section 3.1  in \cite{B}).
\qed

\section{Fluctuations of the zeta zeros on the critical line}

An application of the multidimensional Selberg theorem is to counting zeros of the Riemann zeta functions. Denote by $N(t)$ the number of non-trivial zeros of $\zeta (s) $ on the critical line $\Re s = \frac{1}{2}$ with the imaginary part contained in the interval  $[0,t]$.  Then there is known (see e.g. \cite{Titch}) that
\begin{equation}
\label{24d-1}
N(t)= \frac{t}{2\pi} \log \frac{ t}{2\pi e} + \frac{1}{\pi} \arg \log \zeta \left( \frac{1}{2} + it\right) + \mathcal{O}(\frac{1}{t}).
\end{equation}
If $t_{1} <t_{2}$, let 
\begin{equation}
\label{24d-2}
\Delta (t_{1}, t_{2}) =(N(t_{2})- N(t_{1}) - \left(    \frac{t_{2}}{2\pi} \log \frac{ t}{2\pi e} -\frac{t_{1}}{2\pi} \log \frac{ t_{1}}{2\pi e}\right).
\end{equation}
The quantity $\Delta (t_{1}, t_{2}) $ is  usually  interpreted as the {\it fluctuation of the number of zeta zeros} on the critical line between the heights $\operatorname{Im} s =t_{1}$ and $\operatorname{Im} s = t_{2}$ minus its expected value. 

From the results in the previous section, we can deduce the asymptotic behavior of  the fluctuations of zeta zeros between random points.  

\begin{prop} Let $\Delta $ be given by (\ref{24d-2}) and $U$ by (\ref{u}). Then for every $0<i_{1}< i_{2}$

\begin{equation*}
\frac{1}{\pi  \sqrt{\log \log T}} \Delta \left( UT+ i_{1} T, UT + i_{2}T \right) \to  ^ {(d)}_{T \to \infty} N(0,1)
\end{equation*}
and
\begin{equation*}
d _{W} \left( \frac{1}{\pi \sqrt{\log \log T} } \Delta \left( UT+ i_{1} T, UT + i_{2}T \right) , N(0,1) \right) \leq C \frac{1}{ \sqrt{\log \log T}}.
\end{equation*}

\end{prop}
{\bf Proof: } We proved in Theorem \ref{t6}  that the random vector

$$\left(  \frac{1}{ \sqrt{\frac{1}{2}\log \log T}} (\arg \log \zeta \left( \frac{1}{2} + \mathbf{i}T(U+i) \right) \right) _{i=0,..., d}$$
converges in distribution  to  $ N(0, I_{d+1}) $ with speed less than  $C\frac{1}{\sqrt{\log \log T}}.$ This implies the conclusion. \qed

\vskip0.2cm

In particular, by choosing $i_{1}=0$ and $i_{2}=1$, we will have
\begin{equation}
\label{28d-7}
d_{W}   \left( \frac{1}{\pi \sqrt{\log \log T} } \Delta ( UT, UT+T), N(0,1)\right) \leq C \frac{1}{ \sqrt{\log \log T}}.
\end{equation}

We also have seen that in Theorem \ref{t6} that, if $f_{T} ^ {(i)}$ are as in the statement of Proposition \ref{p4}, then the random vector 

$$\left(  \frac{1}{ \sqrt{\frac{1}{2}\log \log T}} (\arg \log \zeta \left( \frac{1}{2} + \mathbf{i}(UT+ f_{T} ^ {(i)}\right) \right) _{i=0,..., d}$$
converges in distribution to $N(0,\Lambda )$ (the matrix $\Lambda$ has been introduced in Proposition \ref{p4}) at rate $C \frac{1}{\sqrt{\log \log T}}$. Consequently, we have the following result.

\begin{prop}

 For $0<i_{1} <i_{2} $ (so $f_{T} ^ {(i_{1})}< f_{T} ^ {(i_{2})}$) with $f_{T} ^ {(i)}$, $i=0,..,d$ satisfying the assumptions in Proposition \ref{p4},

\begin{equation*}
\frac{1}{\pi \sqrt{\log \log T} } \Delta \left( UT + f_{T} ^ {(i_{1})}, UT + f_{T} ^ {(i_{2})}\right) \to N(0, 1-c_{i_{1},i_{2}}).
\end{equation*}
and
\begin{equation*}
d_{W} \left( \frac{1}{\pi \sqrt{\log \log T} } \Delta \left( UT + f_{T} ^ {(i_{1})}, UT + f_{T} ^ {(i_{2})}\right), N(0, 1-c_{i_{1},i_{2}})\right) \leq C \frac{1}{\sqrt{\log \log T}}.
\end{equation*}
\end{prop}
If we choose $d=1$, $ f_{T}^ {(0)} =0$ and $f_{T} ^ {(1)} = \frac{1}{ (\log T) ^ {\delta } } $ with $0<\delta <1$, then 
\begin{equation*}
 \frac{ \log \vert f_{T} ^ {(1) }- f_{T} ^ {(0)} \vert } {-\log \log T} \to _{T} c_{1,0}:= \delta
\end{equation*}
so

$$ \frac{ \Delta\left( (UT, UT +   \frac{1}{ (\log T) ^ {\delta } }) \right)}{\pi \sqrt{\log \log T}}  \to \sqrt{ 1-\delta } N$$
and
\begin{equation}
\label{28d-8}
d_{W} \left( \frac{ \Delta\left( (UT, UT +   \frac{1}{ (\log T) ^ {\delta } }) \right)}{\pi \sqrt{\log \log T}} ; \sqrt{ 1-\delta } N\right) \leq C \frac{1}{ \sqrt{\log \log T}}.
\end{equation}
Relations (\ref{28d-7}) and (\ref{28d-8}) can be interpreted as follows. The number of zeta zeros on the critical line between the heights $UT$ is $UT+T$ is " close" to $\sqrt{\log \log T}$. With the same approximation error, the number of zeta zeros between the horizontal lines $UT$ and $UT + \frac{1}{ (\log T) ^ {\delta } }$ is approximately $\sqrt{1-\delta }\sqrt{ \log \log T}$ with $\delta \in (0,1)$.

A last consequence concerns the so-called {\it mesoscopic flutuations } of the zeta zeros. 
 
\begin{corollary}
If $K_{T}$ is a deterministic sequence such that $K_{T}>\varepsilon >0$ for every $T>0$ and 
\begin{equation*}
\frac{\log K_{T} }{\log \log T} \to_{T \in \infty} \delta \in [0,1)
\end{equation*}
then the process 
$$\left(\frac{ \Delta (UT+ \frac{\alpha }{K_{T}}), \Delta (UT+ \frac{\beta }{K_{T}})}{ \frac{1}{\pi} \sqrt{ \frac{1}{2} (1-\delta) \log \log T}}, 0\leq \alpha <\beta < \infty\right)$$
converges in the sense of finite dimensional distributions to the centered  Gaussian process $\left( G(\alpha, \beta), 0\leq \alpha <\beta <\infty \right)$ with covariance 
$$ \mathbf{E} G(\alpha , \beta) G( \alpha ' , \beta ' ) =1_{ \left( (\alpha =\alpha ') \mbox{ and } (\beta = \beta ' )\right)}+ \frac{1}{2} 1 _{\left( (\alpha =\alpha ') \mbox{ and } (\beta \not= \beta ' )\right)}+ \frac{1}{2} 1 _{\left( (\alpha \not=\alpha ') \mbox{ and } (\beta = \beta ' )\right)}-\frac{1}{2} 1_{ (\beta = \alpha ' )}.$$

The Wassestein distance associated to this convergence is of order less than $ C \frac{1}{\sqrt{\log \log T}}$.
\end{corollary}

This is interpreted in \cite{B} or \cite{CoDi} as a mesoscopic repulsion of zeros. (Recall that mesoscopic means at a scale between microscopic and macroscopic.) The result shows that the zeta zeros do not affect too much the behavior of $\zeta$  on the critical line.

\section{Appendix}

\subsection{Elements of number theory}

Let $\pi (x)$ be the prime-counting function that gives the number of primes less than or equal to x, for any real number x.  The prime number theorem then states that $\pi (x)$ behaves, when $x$ is large, as  $\frac{x}{\log x}$. As a consequence of this result, certain partial sums of primes can be estimated. We list below some estimates that are needed in our work.  

For every $s$ with $\Re s <1$ we have

\begin{equation}\label{pn1}
\sum _{p\leq x}p ^{-s}\sim \frac{x^{1-s}}{(1-s) \log s}
\end{equation}
while if $s=1$, the sum of the reciprocals of primes diverges as
\begin{equation}\label{pn2}
\sum _{p\leq x}\frac{1}{p}= \log \log x + C + \mathcal{O}\left( \frac{1}{\log x}\right).
\end{equation}
We will also use (see e.g. \cite{Ts})

\begin{equation}\label{pn3}
\sum _{p\leq x}\frac{\log p}{p}\sim   \log x
\end{equation}
and

\begin{equation}\label{pn4}
\sum _{p\leq x}\log p\sim    x.
\end{equation}

\subsection{Basics of the Malliavin calculus}\label{mal}

We present  the elements from the Malliavin calculus  that we  need in the paper.   Consider ${\mathcal{H}}$ a real separable Hilbert space and $(B (\varphi), \varphi\in{\mathcal{H}})$ an isonormal Gaussian process \index{Gaussian process} on a probability space $(\Omega, {\cal{A}}, P)$, which is a centered Gaussian family of random variables such that $\mathbf{E}\left( B(\varphi) B(\psi) \right)  = \langle\varphi, \psi\rangle_{{\mathcal{H}}}$.

We denote by $D$  the Malliavin  derivative operator that acts on smooth functions of the form $F=g(B(\varphi _{1}), \ldots , B(\varphi_{n}))$ ($g$ is a smooth function with compact support and $\varphi_{i} \in {{\cal{H}}}, i=1,...,n$)
\begin{equation*}
DF=\sum_{i=1}^{n}\frac{\partial g}{\partial x_{i}}(B(\varphi _{1}), \ldots , B(\varphi_{n}))\varphi_{i}.
\end{equation*}
It can be checked that the operator $D$ is closable from $\mathcal{S}$ (the space of smooth functionals as above) into $ L^{2}(\Omega; \mathcal{H})$ and it can be extended to the space $\mathbb{D} ^{1,p}$ which is the closure of $\mathcal{S}$ with respect to the norm
\begin{equation*}
\Vert F\Vert _{1,p} ^{p} = \mathbf{E} F ^{p} + \mathbf{E} \Vert DF\Vert _{\mathcal{H}} ^{p}. 
\end{equation*}
By $L$ we will denote the infinitesimal genetaror of the ornstein-Uhlenbeck semigroup and by $L^ {-1}$ its pseudo-inverse. The reader may consult the monographs \cite{N} or \cite{NPbook} for the definition and the properties of this operators. What is need in this paper concerning $L$ and $L^ {-1}$ is only the formula (\ref{nv}).

\end{document}